\documentclass[10pt]{tran-l}
\usepackage{amssymb,amsmath,amsthm,amsfonts,amsopn,url,bbm,enumerate}
\theoremstyle{plain}
\newtheorem{thm}{Theorem}[section]
\newtheorem{prop}[thm]{Proposition}
\newtheorem{lemma}[thm]{Lemma}
\newtheorem{cor}[thm]{Corollary}
\theoremstyle{definition}
\newtheorem{definition}[thm]{Definition}
\newcommand{\field}[1]{\mathbbm{#1}}
\newcommand{\N}{\field{N}}
\newcommand{\Z}{\field{Z}}

\newcommand{\ideal}[1]{\mathfrak{#1}}
\newcommand{\m}{\ideal{m}}

\newcommand{\p}{\ideal{p}}
\newcommand{\q}{\ideal{q}}
\newcommand{\ffunc}[1]{\mathrm{#1}}
\newcommand{\func}[1]{\mathrm{#1} \,}

\newcommand{\Spec}{\func{Spec}}

\newcommand{\Tor}{\ffunc{Tor}}
\newcommand{\Ass}{\func{Ass}}

\newcommand{\coker}{\func{coker}}
\newcommand{\im}{\func{im}}
\newcommand{\bx}[1]{{\mathbf{x}^{[#1]}}}
\newcommand{\bxp}[1]{{\mathbf{x}'^{[#1]}}}
\newcommand{\bfx}{\mathbf{x}}
\newcommand{\bfxp}{\mathbf{x}'}

\newcommand{\pdp}[1]{\ffunc{ph. depth}_{#1}}
\newcommand{\arrow}[1]{\stackrel{#1}{\rightarrow}}
\newcommand{\Assph}{\func{Ass}^{\ffunc{Ph}}}
\newcommand{\GAss}[1]{\bigcup_{e\geq 0} \func{Ass} G^e_R(#1)}
\newcommand{\<}{\left<}
\renewcommand{\>}{\right>}
\title
{Phantom depth and stable phantom exactness.}

\author[Neil Epstein]{Neil M. Epstein}

\address{Department of Mathematics, University of Kansas, Lawrence, Kansas  66045}

\email{epstein@math.ukans.edu}

\thanks{The author was partially supported by the National Science Foundation.}

\keywords{tight closure, phantom depth, phantom homology, rigidity}

\subjclass[2000]{Primary 13A35; Secondary 13C15, 13D25}

\begin{document}
\begin{abstract}
Phantom depth, phantom nonzerodivisors, and phantom exact sequences are analogues of
the non-``phantom'' notions which have been useful in tackling the (very difficult) localization
problem in tight closure theory.  In the present paper, these notions are
developed further and partially reworked.  For instance, although
no analogue of a long exact sequence arises from a short stably phantom exact sequence
of complexes, we provide a method for recovering the kind of information
obtainable from such a long sequence.  Also, we give alternate characterizations
of the notion of phantom depth, including one based on Koszul homology which we use to show that 
with very mild conditions on 
a finitely generated module $M$, any two maximal phantom $M$-regular sequences in an ideal
$I$ have the same length.  In order to do so, we prove a ``Nakayama lemma for tight closure'' which is
of independent interest.  We strengthen the connection of phantom depth
with minheight, we explore several analogues 
of ``associated prime'' in tight closure theory, and we discuss a connection
with the problem of when tight closure commutes with localization.
\end{abstract}

\maketitle

\setcounter{tocdepth}{1}
\tableofcontents

\section{Introduction}\label{sec:intro}
This work concerns the theory of tight closure (in positive characteristic) developed by M. Hochster and C. Huneke in~\cite{HHmain}.  All rings in this work are commutative, Noetherian, local, and of positive prime characteristic $p>0$, and all $R$-modules are finitely generated unless otherwise noted.  For general references on commutative algebra for terms not explained, I recommend the excellent books by Matsumura \cite{Mats} and Bruns \& Herzog \cite{BH}.

One of the many good properties of the tight closure operation is \emph{colon capturing}.  That is, for any system of parameters $x_1, \dotsc, x_d$ in a local ring $(R,\m)$ satisfying mild conditions, $(x_1, \dotsc, x_{i-1}) : x_i \subseteq (x_1, \dotsc, x_{i-1})^*$ for $i = 1, \dotsc, d$.\cite[Theorems 3.1 and 3.1A]{HuTC}  In Mel Hochster's words, tight closure ``captures the failure'' of a system of parameters to be a regular sequence for $R$. Said another way, any system of parameters for $R$ is an $R$-regular sequence ``up to tight closure''.  It is natural to ask the question: For an $R$-module $M$, what sort of sequences in an ideal $I$ are, similarly, $M$-regular ``up to tight closure''?

Aberbach introduces this question in \cite{AbPPD} and assigns the description \emph{phantom $M$-regular} to such a sequence, appropriately defined.  In addition to basic properties of phantom $M$-regular sequences and their companion invariant, \emph{phantom depth}, Aberbach also connects the notion from the title of his paper (\emph{phantom projective dimension}) with phantom depth in an elegant formula.  Namely, he proves a \emph{phantom Auslander-Buchsbaum theorem}, in which, with mild conditions on the ring, the word ``phantom'' is placed in front of every occurrence of the word ``depth'' and the phrase ``projective dimension'' in the statement of Auslander and Buchsbaum's classical theorem.  The phantom projective dimension of $M$ is the length of the shortest \emph{stably phantom acyclic complex} (i.e. a complex which is acyclic up to tight closure and Frobenius powers; see \cite{HHmain} and Section~\ref{sec:pseudolesABBR} of this paper) of finitely-generated projective modules whose zero'th homology module is $M$.  Aberbach, Hochster, and Huneke \cite{AHH} use Aberbach's phantom Auslander-Buchsbaum theorem to obtain some of the best results to date on the difficult localization problem in tight closure theory.

In this work, we investigate phantom depth and stable phantom properties of complexes further.  Nearly all of the previously existing theorems that use phantom depth rely on the phantom Auslander-Buchsbaum theorem, which means that they only apply to modules of finite phantom projective dimension.  In several cases (see Section~\ref{sec:phass}), I have been able to remove this restriction.  Also, I am able to give a positive answer (in Corollary~\ref{cor:flex}) to a question posed by Aberbach: With mild conditions on the ring, all maximal phantom $M$-regular sequences in a proper ideal $I$ have the same length.

My results rest on three
fundamental tools.  First, I develop a technique for treating a short stably phantom exact sequence of complexes in a way analogous to the way one treats a short exact sequence of complexes.  I do not obtain a long exact sequence, or even a long stably phantom exact sequence, but I can mimic the results we would expect if we had such a long sequence.  The minimal amount of this technique necessary for the bulk of the paper is developed in Section~\ref{sec:pseudolesABBR}.  I give a more detailed version, for both tight and Frobenius closures (in Sections~\ref{sec:pseudolesAPX} and \ref{sec:frobdiag}, respectively).  Secondly, I prove tight closure versions of the Nakayama lemma (in Section~\ref{sec:Nak*}).  Thirdly, I define ``ghost $M$-regular sequences'', which coincide with phantom $M$-regular sequences in all cases of interest (see Proposition~\ref{prop:gh=ph}) but are often easier to use than the latter.  With these tools, I obtain characterizations of phantom depth (Theorem~\ref{thm:flex}) and phantom $M$-regular sequences (Theorem~\ref{thm:ghseqK}) in terms of stable phantomness, rather than vanishing, of Koszul homology.
In Section~\ref{sec:rtexact}, I give an analogue of right exactness and show that it interfaces with 
phantom regular sequences in the expected way.

\subsection{Tight closure background} In order to define tight closure, we first need to explain about Frobenius powers:
For a finitely generated module $M$ over such a ring $R$ with maximal ideal $\m$, and for any nonnegative integer $e \geq 0$, we denote the following concept by the symbol $F^e(M)$: Let ${}^e R$ be the additive abelian group $R$ whose $R$-$R$ bimodule structure is given by $a \cdot z \cdot b = ab^q z$, where $a, b \in R$ and $z \in {}^e R$.  Then $F^e(M)$ is the (left) $R$-module ${}^e R \otimes_R M$.  Clearly this makes $F^e$ a right-exact functor on the category of $R$-modules.
Now, let \[
R^n \arrow{\phi} R^m \arrow{\pi} M \rightarrow 0
\]
be a finite free presentation of $M$.  Let $q = p^e$.  If we fix bases for the free modules $R^n$ and $R^m$, we
get a corresponding representation of $\phi$ as a matrix $A = (a_{ij})$ with entries in $\m$.  Let $A^{[q]}$ be the matrix whose $(i,j)$'th entry is $a_{ij}^q$, for all $1 \leq i \leq m$, $1 \leq j \leq n$.  Then $F^e(\phi): R^n \rightarrow R^m$ is the map of free modules represented by the matrix $A^{[q]}$, so that, since $F^e$ is right-exact, $F^e(M) = \coker F^e(\phi)$.  Moreover, for any $z \in M$, there exists $y \in R^m$ such that $\pi(y) = z$.  Let $\pi': R^m \rightarrow F^e(M) = \coker F^e(\phi)$ be the natural surjection from $R^m$ induced by $F^e(\phi)$.  Then set $z^q := \pi'(y)$.  The notation would be more precise if we were to write $z^q_M$ instead of $z^q$.  However, we shall use the convention that for any element $z$, $z^q$ will denote $z^q_M$, where $M$ is the largest module of those mentioned that contain $z$ as an element.  Other than that, all of the foregoing is functorial and independent of all choices.  If $N$ is a submodule of $M$, let $i: N \rightarrow M$ be the inclusion.  Then $N^{[q]}_M$, the ``$q$'th Frobenius power of $N$ in $M$,'' is defined to be the image of the map $F^e(i): F^e(N) \rightarrow F^e(M)$.

As for tight closure: If $N \subseteq M$ are finitely-generated $R$-modules, then
we say that an element $z \in M$ is in the \emph{tight closure of $N$ with respect to $M$}, and we write $z \in N^*_M$,
if there exists some $c \in R$, which is not in any minimal prime of $R$, such that for all sufficiently large powers
$q$ of $p$, $c z^q \in N^{[q]}_M$.  For any $e \geq 0$, we let $G^e(M) = F^e(M) / 0^*_{F^e(M)}$, and we call it the \emph{$q$'th reduced Frobenius power of $M$}.  If there is some $c \in R$, not in any minimal prime of $R$, and some power $q_0$ of $p$ such that for \emph{all} finitely-generated $R$-modules $M$ and \emph{all} submodules $N \subseteq M$, $z \in N^*_M$ if and only if $c z^q \in N^{[q]}_M$ for all powers $q$ of $p$ such that $q \geq q_0$, then we say that $c$ is a \emph{$q_0$-weak test element} for $R$.  If $c$ is a $q_0$-weak test element for $R_\p$ for every $\p \in \Spec R$, then it is called a \emph{locally stable} $q_0$-weak test element for $R$.  If there exists some power $q_0$ of $p$ for which $c$ is a (locally stable) $q_0$-weak test element, we say that $c$ is a (locally stable) weak test element.  One of the loveliest accomplishments of tight closure theory is the theorem of Hochster and Huneke \cite[Theorem 7.32]{HHbase} which says among other things that if $R$ is excellent and local (which is not a very restrictive condition; for instance any complete local ring is excellent, as is any algebra essentially of finite type over a field), then it has a locally stable weak test element.  Unless otherwise noted, all rings in this paper will contain a $q_0$-weak test element $c$.

\section{Coping without authentic long exact sequences}\label{sec:pseudolesABBR}
The fact that a short exact sequence of complexes
leads to a long exact sequence is
a central tool in the homological theory of modules over
a ring.  When investigating phantom regular sequences
(as we begin to do in Section~\ref{sec:invar})
we obtain sequences which are right-exact and ``stably phantom exact'' on the left.  In this section we derive just the information (the content of Proposition~\ref{prop:pseudolesABBR}) that we need from such a sequence to prove Lemma~\ref{lem:cohoind}, a key step in the proof of Theorem~\ref{thm:flex}, a main result of this paper.  This information is similar to some of what we would obtain
from a true long exact sequence.  In Section~\ref{sec:pseudolesAPX}, we present a more
fully developed theory.

\noindent \textit{Notation and conventions:} Throughout this paper, we shall use the following convention linking the letters $e$ and $q$: any variation of the letter $q$ will be understood to be $p$ to the power of $e$ with the corresponding variation.  For example, $q_0 = p^{e_0}$, $q_1 = p^{e_1}$, $q' = p^{e'}$, and $q'' = p^{e''}$.

In a notational abbreviation which we shall use throughout the paper, given $q_0$ and a fixed $q_0$-weak test element $c$, recursively define $c_n$ for integers $n \geq -1$ by the rules: \[
c_{-1} = 1
\] and \[
c_{n+1} = c \cdot c_n^{\ q_0},
\]
In particular, $c_0 = c$, $c_1 = c^{q_0+1}$,
$c_2 = c^{q_0^2 + q_0 + 1}$, etc.  This is useful because whenever
$A \subseteq B$ are finitely generated $R$-modules and $c_n z \in A^*_B$, it follows that 
$c_{n+1} z^{q_0} = c (c_n z)^{q_0} \in A^{[q_0]}_B$.

If $M.$ is a complex of $R$-modules, then the $i$'th differential is denoted $d_i^M$, the $i$'th module
of cycles is $Z_i(M.) := \ker d_i^M$ and the $i$'th boundary module is $B_i(M.) := \im d_{i+1}^M$.

In this section and in Section~\ref{sec:pseudolesAPX}, to simplify notation, we will abuse it by assuming the
reader can keep track of the
Frobenius powers on the maps.  In particular, if $g: X \rightarrow Y$ is a map of $R$-modules,
$g$ will denote $F^e(g)$ 
for some $e$, and if $M.$ is a complex of $R$-modules, 
$d_i^M$ will denote some Frobenius power of the $i$'th differential of the complex $M$.  Also in
this section, we will mix homology with Frobenius powers:
If $\alpha.: L. \rightarrow M.$ is a map of complexes of $R$-modules and $i \in \Z$, then $H_i(F^e(\alpha.))$
denotes the map $H_i(F^e(L.)) \rightarrow H_i(F^e(M.))$ induced by the composition of the 
$i$'th homology functor $H_i$ with the Frobenius functor $F^e$ acting on the map $\alpha.$ of complexes.  An element of $H_i(F^e(\alpha.))$ will
be denoted in brackets, e.g.: $[x]$.  Combining all these conventions together, if $x \in \ker F^e(d_i^L)$,
then we write $[x] \in H_i(F^e(L.))$, and the symbols $H_i(F^e(\alpha.))([x])$ and $[\alpha_i(x)]$ denote
the same element of $H_i(F^e(M.))$.

\begin{definition}\cite{HHmain}
If $M.$ is a complex of $R$-modules and $i \in \Z$,  we say that \emph{$M.$
has stably phantom homology at the $i$'th spot} (or \emph{at $i$}) if \[
Z_i(F^e(M.)) \subseteq \left( B_i(F^e(M.)) \right)^*_{F^e(M_i)}.
\]
for all $e \geq 0$.  We say that a complex is \emph{stably phantom exact}
(resp. \emph{stably phantom acyclic})
if it has stably phantom homology at every spot (resp. at every spot except
0).

We say that an element $[x]$ of $H_i(M.)$ is \emph{phantom}
if for all $e \gg 0$, $c x^q \in \im F^e(d^{i+1}_M)$ (i.e. $[c x^q] = 0$ in 
$H_i(F^e(M_.))$).
\end{definition}

At this point I collect together two easily proved Facts which will be used in the
sequel without comment.

\begin{enumerate}
\item [Fact 1:] For any map $\beta : L \rightarrow N$ of finitely
generated $R$-modules and any $e$, $\im F^{e}(\beta) = (\im
\beta)^{[q]}_{N}$.

\item [Fact 2:] For map $\beta : L \rightarrow N$ of finitely
generated $R$-modules and any $e$, $(\ker \beta)^{[q]}_{L} \subseteq
\ker F^{e}(\beta)$.
\end{enumerate}

Next, we have some characterizations of stably phantom homology:

\begin{prop}\label{prop:killphantom}
Let $R$ be a Noetherian ring of prime characteristic $p>0$ which contains a $q_0$-weak test element
$c$, let $L.$ be a complex of finitely generated $R$-modules.  For any $i \in \Z$, the following
are equivalent: \begin{enumerate}
\item[(a)]\label{it:phantom} $L.$ has stably phantom homology at $i$.
\item[(b)]\label{it:0*coker} $H_i(F^e(L.)) \subseteq 0^*_{\coker d_{i+1}^{F^e(L.)}}$ for all $e \geq 0.$
\item[(c)]\label{it:ckill} For all $e \geq 0$, $c$ kills $(H_i(F^e(L.)))^{[q_0]}_{\coker
 F^e(d_{i+1}^L)}$.
\end{enumerate}
\end{prop}

\begin{proof}
The equivalence of (a) and (b) is easy, so we prove the equivalence
of (a) and (c).

If $L.$ has stably phantom homology at $i$, then for any $e$ and any $[z] \in H_i(F^e(L.))$, we have that
$c z^{q_0} \in \im 
F^{e+e_0}(d_{i+1}^L)$, which means that $c [z]^{q_0}_{\coker F^e(d_{i+1}^L)} = [c z^{q_0}] = 0$.

Conversely, suppose (c) is true.  Take any $[x] \in H_i(F^e(L.))$.  Then 
for all powers $q'$ of $p$, $[x^{q'}] \in H_i(F^{e+e'}(L.))$, so that by (c),
$[c x^{q'q_0}] = 0$.  Since this holds for all large $q'$, it follows that $[x]$ is a
phantom element of $H_i(F^e(L.))$.
\end{proof}

The following is the main result of Section~\ref{sec:pseudolesABBR}.

\begin{prop}\label{prop:pseudolesABBR}
Let $R$ be a Noetherian ring of prime characteristic $p>0$ which contains a $q_0$-weak test element
$c$, and let \[
0 \rightarrow L. \stackrel{\alpha.}{\rightarrow} M. \stackrel{\beta.}{\rightarrow} N. \rightarrow 0
\]
be a right-exact sequence of complexes of finitely generated $R$-modules, such that
for each degree $i$ and every integer $e \geq 0$, $\ker F^e(\alpha_i) \subseteq
0^*_{F^e(L_i)}$. (i.e. the sequence is both right exact and stably phantom exact)

Then for any $i \in \Z$: 
\begin{enumerate}
\item\label{it:Nph} $N.$ has stably phantom homology at $i$ if and only if the
following two conditions hold for all $e \geq 0$:
\begin{enumerate}
  \item Any element of $H_{i-1}(F^e(L.))$ which $H_{i-1}(F^e(\alpha.))$ maps
  to a phantom element of $H_{i-1}(F^e(M.))$ is itself a phantom element of 
  $H_{i-1}(F^e(L.))$, and
  \item For any $[y] \in H_i(F^e(M.))$ and any integer $e' \geq e_0$,
  $[c y^{q'}] \in \im H_i(F^{e+e'})(\alpha.)$.
\end{enumerate}
\item\label{it:LMph} If $M.$ has stably phantom homology at $i$ and $L.$ has stably phantom homology at $i-1$,
then $N.$ has stably phantom homology at $i$.
\end{enumerate}
\end{prop}

\noindent \emph{Note:}  Part (1) of the above Proposition is an analogue of the fact that in a 
five-term exact sequence, the middle term is zero if and 
only if the first map is surjective and the last map is injective. Part (2) is an analogue
of the fact that in a three-term exact sequence, if the outer two terms vanish, so does the
middle term.

Also note that the result of applying the functor $F^e$ to the sequence given in the Proposition
preserves the hypotheses.  For instance, the new sequence of complexes is right-exact because
the old one was and $F^e$ is a right-exact functor.

These proofs arose as generalizations of proofs of their non-phantom analogues.

\begin{proof}[Proof (\ref{it:Nph})]
First suppose that $N.$ has stably phantom homology at $i$.

To prove condition (a), let $[x] \in H_{i-1}(F^e(L.))$ be such that 
$[\alpha_{i-1}(x)]$ is a phantom element of $H_{i-1}(F^e(M.))$.  Then
for any $e' \geq e_0$, there exists $y \in F^{e+e'}(M_i)$ such that
$\alpha_{i-1}(c x^{q'}) = c (\alpha_{i-1}(x))^{q'} = d_i^M(y)$.  Then
\[
d_i^N(\beta_i(y)) = \beta_{i-1}(d_i^M(y)) = \beta_{i-1}(\alpha_{i-1}(c x^{q'})) = 0,
\]
so $[\beta_i(y)] \in H_i(F^{e+e'}(N))$.  By stable phantomness of $N.$ at $i$ and
surjectivity of $\beta_{i+1}$,
there exists $v \in F^{e+e'+e_0}(M_{i+1})$ with \[
c \beta_i(y)^{q_0} = \beta_i(c y^{q_0}) = d_{i+1}^N(\beta_{i+1}(v)) = \beta_i(d_{i+1}^M(v)).
\]
Thus, $c y^{q_0} - d_{i+1}^M(v) \in \ker \beta_i = \im \alpha_i$, so that 
there is some $u \in F^{e+e'+e_0}(L_i)$ with $c y^{q_0} = d_{i+1}^M(v) + \alpha_i(u)$.
We have: \begin{align*}
\alpha_{i-1}(c_1 x^{q'q_0}) &= c (\alpha_{i-1}(c x^{q'}))^{q_0} = 
c (d_i^M(y))^{q_0} \\
&= d_i^M(c y^{q_0}) = d_i^M(\alpha_i(u)) = \alpha_{i-1}(d_i^L(u)).
\end{align*}
Thus, $c_1 x^{q' q_0} - d_i^L(u) \in \ker \alpha_{i-1} \subseteq 0^*_{L_{i-1}}$,
so that \[
c_2 x^{q' q_0^2} = d_i^L(c u^{q_0}) \in (\im d_i^L)^{[q' q_0^2]}_{F^e(L_{i-1})}.
\]
Since this holds for any $e' \geq e_0$, it follows that $x \in (\im d_i^L)^*_{F^e(L_{i-1})}$,
whence $[x]$ is a phantom element of $H_{i-1}(F^e(L.))$.  

As for condition (b), note that for any $[y] \in H_i(F^e(M.))$, $[\beta_i(y)]$ is 
phantom by hypothesis, so for any $e' \geq e_0$, since $\beta_{i+1}$ is surjective
there exists $v \in F^{e+e'}(N_{i+1})$ with \[
\beta_i(c y^{q'}) = d_{i+1}^N(\beta_{i+1}(v)) = \beta_i(d_{i+1}^M(v)),
\]
which implies that $c y^{q'} - d_{i+1}^M(v) \in \ker \beta_i = \im \alpha_i$.  That is, 
there exists $u \in F^{e+e'}(L_i)$ with $c y^{q'} = d_{i+1}^M(v) + \alpha_i(u)$,
which proves condition (b).

Conversely, suppose that conditions (a) and (b) hold for all $e \geq 0$, and
let $[z] \in H_i(F^e(N.))$.  Then $z = \beta_i(y)$ for some $y \in F^e(M_i)$,
since $\beta_i$ is surjective.  Moreover, $0 = d_i^N(z) = d_i^N(\beta_i(y)) 
= \beta_{i-1}(d_i^M(y))$, so that $d_i^M(y) \in \ker \beta_{i-1} = \im \alpha_{i-1}$.
Hence, there exists $x \in F^e(L_{i-1})$ with $\alpha_{i-1}(x) = d_i^M(y)$.  We 
have \[
\alpha_{i-2}(d_{i-1}^L(x)) = d_{i-1}^M(\alpha_{i-1}(x)) = d_{i-1}^M(d_i^M(y)) = 0,
\]
so that $d_{i-1}^L(x) \in \ker \alpha_{i-2} \subseteq 0^*_{F^e(L_{i-2})}$.  Thus
for any $e' \geq e_0$, we have $d_{i-1}^L(c x^{q'}) = c (d_{i-1}^L(x))^{q'} = 0$, so that
$[c x^{q'}] \in H_{i-1}(F^{e+e'}(L.))$.  Moreover, $H_{i-1}(F^{e+e'}(\alpha))([c x^{q'}]) =
[\alpha_{i-1}(c x^{q'})] = [d_i^M(c y^{q'})] = 0$, so that by condition (a), 
$[c x^{q'}]$ is phantom, which implise that the existence of some $t \in F^{e+e'+e_0}(L_i)$
with $c_1 x^{q' q_0} = d_i^L(t)$.  Now, \[
d_i^M (\alpha_i(t)) = \alpha_{i-1}(d_i^L(t)) = \alpha_{i-1}(c_1 x^{q' q_0}) = d_i^M(c_1 y^{q' q_0}),
\]
so that $[c_1 y^{q' q_0} - \alpha_i(t)] \in H_i(F^{e+e'+e_0}(M.))$, from which
condition (b) implies that \[
[c_2 y^{q' q_0^2} - \alpha_i(c t^{q_0})] = \left[c \left( c_1 y^{q' q_0} - \alpha_i(t) 
 \right)^{q_0}\right] \in \im H_i(F^{e+e'+2e_0}(\alpha)).
\]
That is, there exist $u \in F^{e+e'+2e_0}(L_i)$ and $v \in F^{e+e'+2e_0}(M_{i+1})$
such that $c_2 y^{q' q_0^2} - \alpha_i(c t^{q_0}) =  \alpha_i(u) + d_{i+1}^M(v)$.
Hence, \[
c_2 z^{q' q_0^2} = \beta_i(c_2 y^{q' q_0^2}) = \beta_i(d_{i+1}^M(v))
 = d_{i+1}^N(\beta_{i+1}(v)).
\]
Therefore $[c_2 z^{q' q_0^2}] = 0$, which means that $[z]$ is phantom, proving that
$N.$ has stably phantom homology at $i$.
\end{proof}

\begin{proof}[Proof (\ref{it:LMph}).]
Let $e \geq 0$ be arbitrary.  Since $L.$ has stably phantom homology at $i-1$,
condition (\ref{it:Nph}a) holds.  For any $[y] \in H_i(F^e(M.))$, $[c y^{q'}] = 0
= H_i(F^{e+e'}(\alpha.))(0)$ for all $e' \geq e_0$, which proves condition (\ref{it:Nph}b).
Thus by part (\ref{it:Nph}), $N.$ has stably phantom homology at $i$.
\end{proof}

\section{Nakayama lemmas for tight closure}\label{sec:Nak*}
We prove here tight closure versions of the Nakayama lemma, which we use
both in this paper (Sections~\ref{sec:invar} and \ref{sec:phseqK}) and also in
a very different way in \cite{nme*spread}.

\begin{prop}[Nakayama lemma for tight closure, nuts-and-bolts version]\label{prop:Nak*a}
Let $(R,\m)$ be a Noetherian local ring of characteristic $p>0$ possessing
a $q_0$-weak test element $c$.  Let $L$ be a submodule of $M$, and for each
integer $e \geq 0$ let $N_e$ be a submodule of $F^e(M)$ which contains 
$L^{[q]}_M$.  (In particular, this last condition is satisfied as long as $L \subseteq N_0$
and $(N_0)^{[q]}_M \subseteq N_e$ for all $e \geq 0$.)  Suppose in addition that 
for all $e \geq 0$ and $e' \geq e_0$, we have
\begin{equation}\label{eq:containm}
c (N_e)^{[q']}_{F^e(M)} \subseteq L^{[q q']}_M + \m^{[q q']} N_{e+e'}.
\end{equation}
Then $N_e \subseteq (L^{[q]}_M)^*_{F^e(M)}$ for all $e \geq 0$.
\end{prop}

\begin{proof}
As a first step, 
we show the following by induction on $r$:

\noindent \textbf{Claim: } For all integers $r \geq 0$, $e \geq 0$ and $e' \geq e_0$,
\begin{equation}\label{eq:containmr}
c_r (N_e)^{[q' q_0^r]}_{F^e(M)} \subseteq L^{[q q' q_0^r]}_M 
+ \left(\m^{[q q' q_0^r]}\right)^{r+1} N_{e+e'+re_0}.
\end{equation}

\begin{proof}[Proof of Claim]
The case $r=0$ is true by hypothesis.  So let $r>0$ and assume that (\ref{eq:containmr})
is true for $r-1$.  Then \begin{equation}\label{eq:1}
c_{r-1} (N_e)^{[q' q_0^{r-1}]}_{F^e(M)} \subseteq L^{[q q' q_0^{r-1}]}_M 
+ \left(\m^{[q q' q_0^{r-1}]}\right)^r N_{e+e'+(r-1)e_0}
\end{equation}
by inductive hypothesis, and \begin{equation}\label{eq:2}
c (N_{e+e'+(r-1)e_0})^{[q_0]}_{F^{e+e'+(r-1)e_0}(M)} \subseteq L^{[q q' q_0^r]}_M
+ \m^{[q q' q_0^r]} N_{e + e' + r e_0}
\end{equation}
by replacing $e$ by $e+e'+(r-1)e_0$ and $e'$ by $e_0$ in (\ref{eq:containm}).

Now apply the operator $c (-)^{[q_0]}_{F^{e+e'+(r-1)e_0}(M)}$ to both sides
of (\ref{eq:1}) to get: \begin{align*}
c_r (N_e)^{[q' q_0^r]}_{F^e(M)} &= c \left(c_{r-1} 
 (N_e)^{[q' q_0^{r-1}]}_{F^e(M)} \right)^{[q_0]}_{F^{e+e'+(r-1)e_0}(M)} \\
&\subseteq c \left(L^{[q q' q_0^{r-1}]}_M 
 + \left(\m^{[q q' q_0^{r-1}]}\right)^r N_{e+e'+(r-1)e_0}
 \right)^{[q_0]}_{F^{e+e'+(r-1)e_0}(M)} \\
&= c L^{[q q' q_0^r]}_M + \left(\m^{[q q' q_0^r]}\right)^r 
 \left( c (N_{e+e'+(r-1)e_0})^{[q_0]}_{F^{e+e'+(r-1)e_0}(M)} \right) \\
&\subseteq c L^{[q q' q_0^r]}_M + \left(\m^{[q q' q_0^r]}\right)^r 
 \left(L^{[q q' q_0^r]}_M + \m^{[q q' q_0^r]} N_{e + e' + r e_0} \right) \\
&= \left(c + \left(\m^{[q q' q_0^r]}\right)^r \right) L^{[q q' q_0^r]}_M 
 + \left(\m^{[q q' q_0^r]}\right)^r \m^{[q q' q_0^r]} N_{e + e' + r e_0} \\
&\subseteq L^{[q q' q_0^r]}_M + \left( \m^{[q q' q_0^r]} \right)^{r+1}
 N_{e + e' + r e_0},
\end{align*}
proving the Claim.  The containment in the second
line follows from (\ref{eq:1}), while that in 
fourth line follows from (\ref{eq:2}).
\end{proof}

Fixing $r$ and $e$ and letting $e'$ vary, set $e'' = e' + r e_0$ to simplify the 
notation.  Then (\ref{eq:containmr}) says that \begin{align*}
c_r (N_e)^{[q'']}_{F^e(M)} &\subseteq L^{[q q'']}_M + 
 (\m^{[q q'']})^{r+1} N_{e+e''} \\
&\subseteq L^{[q q'']}_M + (\m^{[q q'']})^{r+1} F^{e+e''}(M) \\
&= \left(L^{[q]}_M + (\m^{[q]})^{r+1} F^e(M) \right)^{[q'']}_{F^e(M)}.
\end{align*}
Since this holds for all $e'' \geq (r+1)e_0$, it follows that
\[
N_e \subseteq \left( L^{[q]}_M + (\m^{[q]})^{r+1} F^e(M) \right)^*_{F^e(M)}.
\]

Thus, for any $e' \geq e_0$, we 
have \begin{align*}
c (N_e)^{[q']}_{F^e(M)} &\subseteq \bigcap_{r \geq 0} 
 \left( \left(L^{[q]}_M\right)^{[q']}_{F^e(M)} + (\m^{[q q']})^{r+1} F^{e+e'}(M) \right) \\
&= \left(L^{[q]}_M\right)^{[q']}_{F^e(M)} 
 + \Bigg(\bigcap_{r \geq 0} (\m^{[q q']})^{r+1} \Bigg) F^{e+e'}(M) \\
&= \left(L^{[q]}_M \right)^{[q']}_{F^e(M)},
\end{align*}
where the last equality is a consequence of Krull's intersection theorem.
Hence, $N_e \subseteq (L^{[q]}_M)^*_{F^e(M)}$.
\end{proof}

For aesthetic reasons, and because it is tantamount to the version used
in \cite{nme*spread} to show the existence of ``minimal $*$-reductions'', we include the following corollary here:

\begin{cor}[Nakayama lemma for tight closure, generic version]
Let $(R,\m)$ be a Noetherian local ring of characteristic $p>0$ possessing
a $q_0$-weak test element $c$.  Let
$L \subseteq N \subseteq M$ be finitely generated $R$-modules such that
$N \subseteq (L + \m N)^*_M$.  Then $N \subseteq L^*_M$.
\end{cor}

\begin{proof}
Set $N_e = N^{[q]}_M$ for all $e \geq 0$, and then apply
Proposition~\ref{prop:Nak*a}.
\end{proof}

\newtheorem{av}[thm]{A note on avoidance}
\section{Phantom depth, ghost depth, and lengths of phantom regular sequences}\label{sec:invar}
In this section I introduce Aberbach's phantom analogues from \cite{AbPPD} of depth, zerodivisors, and 
regular sequences.  Then I give alternate characterizations (Propostion~\ref{prop:gh=ph}) of these concepts, and use these
to give yet another characterization (Theorem~\ref{thm:flex}) of phantom depth in terms of stable phantomness of Koszul homology.
From this characterization, we get in Corollary~\ref{cor:flex} a positive answer to a question
Aberbach asked in his paper.

Recall \cite[3.2.1-3.2.2]{AbPPD}:
\begin{definition}\label{def:pdepth}
Let $R$ be a Noetherian ring of characteristic $p>0$ and $M$ a
finitely generated $R$-module.  Then we say an element $x \in R$
is \emph{phantom $M$-regular} (or a \emph{phantom nonzerodivisor
of $M$}) if $x M \neq M$ and $0 :_{F^e(M)} x^t \subseteq 0^*_{F^e(M)}$
for all $e \geq 0$ and all $t \geq 1$.

A \emph{phantom zerodivisor of $M$} is an element of $R$ which
is not phantom $M$-regular.

A sequence $\mathbf{x} = x_1, \dotsc, x_n$ of elements of $R$ is
a \emph{phantom $M$-regular sequence} if $\mathbf{x} M \neq M$ and
if for all $0 \leq i < n$, all $i$-tuples $(u_1, \dotsc, u_i)$ of
positive integers, and all integers $t \geq 1$, $x_{i+1}^t$ is
phantom $(M / (x_1^{u_1}, \dotsc, x_i^{u_i})M)$-regular.
\end{definition}

Recall also that if $\p \in \Spec R$ and $\mathbf{x} = x_1, \dotsc, x_n$ 
is a phantom $M$-regular sequence in $\p$, then $x_1 / 1, \dotsc, x_n / 1$ 
is a phantom $M_\p$ regular sequence.

The following definition appears at first to be strictly weaker than 
Aberbach's definition.  However, under very mild conditions on the ring
(c.f. Proposition~\ref{prop:gh=ph}), the
two notions agree; I know of no instance where they disagree.

\begin{definition}\label{def:gdepth}
Let $R$ be a Noetherian ring of characteristic $p>0$ and $M$ a
finitely generated $R$-module.  Then we say an element $x \in R$ is
\emph{ghost $M$-regular} if $xM \not = M$ and $0:_{F^{e}(M)} x^{q}
\subseteq 0^{*}_{F^{e}(M)}$ for all $e \geq 0$.

A \emph{ghost zerodivisor of $M$} is an element $x \in R$ which is
not ghost $M$-regular.

A sequence $\mathbf{x} = x_{1}, \dotsc , x_{n}$ of elements of $R$ is
a \emph{ghost $M$-regular sequence} if $\mathbf{x}M \not =M$ and
$x_{i+1}$ is ghost $(M / (x_{1}, \dotsc , x_{i})M)$-regular for $0
\leq i < n$.
\end{definition}

Note that an element $x \in R$ with $xM \neq M$ is
ghost $M$-regular if and only if the sequence $0 \rightarrow M 
\stackrel{x}{\rightarrow} M$ is stably phantom exact.  The analogy with
$M$-regular elements is clear.

The first thing to note about ghost $M$-regular sequences is that those
which sit in the maximal ideal of a local ring with a weak test element are permutable:

\begin{prop}\label{prop:permute}
Let $(R,\m)$ be a Noetherian local ring of characteristic $p>0$ 
with a $q_0$-weak test element $c$, and let $M$ be a finitely-generated
$R$-module.  Then for any ghost $M$-regular sequence $\bfx$ contained in $\m$, 
any permutation of $\bfx$ is also a ghost $M$-regular sequence.
\end{prop}
\begin{proof}
It is clear from the definition that for any $1 \leq i < n$, $x_1, \dotsc, 
x_n$ is a ghost $M$-regular sequence if and only if $x_1, \dotsc, x_i$ 
is a ghost $M$-regular sequence and $x_{i+1}, \dotsc, x_n$ is a ghost 
$\left(M / (x_1, \dotsc, x_i)M\right)$-regular sequence.  Hence, since any
permutation is a composition of transpositions of adjacent elements, it 
suffices to show the result for ghost $M$-regular sequences of length two.
That is, we assume that $x, y$ is a ghost $M$-regular sequence, and we need
to show that $y,x$ is a ghost $M$-regular sequence.

First we show that $x$ is ghost $(M / yM)$-regular:  Let $e \geq 0$, 
and let $\bar{z} \in F^e(M / yM) = F^e(M) / 
y^q F^e(M)$ such that
$x^q \bar{z} = \bar{0}$.  That is, $z \in F^e(M)$, and \begin{equation}\label{eq:10*}
x^q z = y^q w
\end{equation}
for some $w \in F^e(M)$.  Hence $y^q w \in x^q F^e(M)$, so since $y$ is
ghost $(M / xM)$-regular, $w \in (x^q F^e(M))^*_{F^e(M)}$, which implies that
for any $e' \geq e_0$, there exists $v \in F^{e+e'}(M)$ with \begin{equation}\label{eq:10**}
c w^{q'} = x^{q q'} v.
\end{equation}
Combining (\ref{eq:10*}) and (\ref{eq:10**}), we get \[
x^{q q'} y^{q q'} v = c y^{q q'} w^{q'} = c \left(y^q w\right)^{q'} 
= c (x^q z)^{q'} = c x^{q q'} z^{q'}.
\]
That is, $x^{q q'} (c z^{q'} - y^{q q'} v) = 0 \in F^{e+e'}(M)$.
By ghost $M$-regularity of $x$, then, \[
c_1 z^{q' q_0} - c y^{q q' q_0} v^{q_0} = c (c z^{q'} - y^{q q'} v)^{q_0} = 0.
\]
Thus, $c_1 z^{q' q_0} \in (y^q F^e(M))^{[q' q_0]}_{F^e(M)}$.  Since $q'$ may be 
arbitrarily large, it follows that $z \in (y^q F^e(M))^*_{F^e(M)}$, which proves
that $x$ is ghost $(M / y M)$-regular.

Now we show that $y$ is ghost $M$-regular.  Let $e \geq 0$ and $z \in F^e(M)$
such that $y^q z = 0$.  Then $y^q z \in x^q F^e(M)$, so that by ghost
$(M / x M)$-regularity of $y$, $z \in \left(x^q F^e(M)\right)^*_{F^e(M)}$.
Hence for any $e' \geq e_0$, there exists $w \in F^{e+e'}(M)$ with
$c z^{q'} = x^{q q'} w$.  We have \[
x^{q q'} y^{q q'} w = c y^{q q'} z^{q'} = c (y^q z)^{q'} = c 0^{q'} = 0,
\]
so that by ghost $M$-regularity of $x$, $y^{q q'} w \in 0^*_{F^{e+e'}(M)}$.

Thus, $c y^{q q' q_0} w^{q_0} = 0$, so that $c w^{q_0} \in 0 :_{F^{e+e'+e_0}(M)} y^{q q' q_0}$. 
Hence, \begin{align*}
c_1 z^{q' q_0} &= c (c z^{q'})^{q_0} = c (x^{q q'} w)^{q_0} = x^{q q' q_0} (c w^{q_0}) \\
&\in x^{q q' q_0} (0 :_{F^{e+e'+e_0}(M)} y^{q q' q_0}).
\end{align*}
Since $z$ was an arbitrary element of $(0 :_{F^e(M)} y^q)$ and since $x \in \m$,
we have that \[
c_1 (0 :_{F^e(M)} y^q)^{[q'']}_{F^e(M)} \subseteq \m^{[q q'']} (0 :_{F^{e+e''}(M)} y^{q q''})
\]
for all $e'' \geq 2e_0$.  Since $c_1$ is a $q_0$-weak test element, setting  
$N_e = 0 :_{F^e(M)} y^q$ and $L = 0$ in Proposition~\ref{prop:Nak*a} shows that \[
0 :_{F^e(M)} y^q \subseteq 0^*_{F^e(M)}
\]
for all $e \geq 0$.  In other words, $y$ is a ghost $M$-regular element.
\end{proof}

Next, we show that in all cases of interest here, the notions of phantom
and ghost $M$-regular sequences coincide, and also that in these cases an element's
phantom $M$-regularity is characterized by $G^e(M)$-regularity for (almost)
all $e$.

\begin{prop}\label{prop:gh=ph}
Let $(R,\m)$ be a Noetherian local ring of characteristic $p>0$ which
contains a $q_0$-weak test element $c$, let
$M$ be a finitely generated $R$-module, and let $\bfx = x_1, \dotsc, x_n$ be
a sequence of elements of $\m$.  Then the following are equivalent:
\begin{enumerate}
\item[(a)] $M \neq \bfx M$, and for all $0 \leq i < n$, all $(i+1)$-tuples
$(t_1, \dotsc, t_{i+1})$ of positive integers, and all integers $e \geq 0$,
\begin{equation}\label{eq:pdepth+def}
0 :_{F^e(M) / (x_1^{t_1}, \dotsc, x_i^{t_i}) F^e(M)} x_{i+1}^{t_{i+1}}
\subseteq 0^*_{F^e(M) / (x_1^{t_1}, \dotsc, x_i^{t_i}) F^e(M)}.
\end{equation}
\item[(b)] $\bfx$ is a phantom $M$-regular sequence.
\item[(c)] $\bfx$ is a ghost $M$-regular sequence.
\item[(d)] $M \neq \bfx M$, and for all $0 \leq i < n$ and all $e \geq 0$,
$x_{i+1}^q$ is $G^e(M/(x_1, \dotsc, x_i) M)$-regular.
\item[(e)] $M \neq \bfx M$, and for all $0 \leq i < n$ and all $e \gg 0$, 
$x_{i+1}^q$ is $G^e(M/(x_1, \dotsc, x_i) M)$-regular.
\end{enumerate}
\end{prop}

\begin{proof}
We will show that (a) $\Rightarrow$ (b) $\Rightarrow$ (c) $\Rightarrow$ (d) $\Rightarrow$ (e) 
$\Rightarrow$ (c) $\Rightarrow$ (a).

It is easy to see that (a) $\Rightarrow$ (b) $\Rightarrow$ (c).  For suppose
$\bfx$ satisfies (a).  Then for any $i$-tuple $(u_1, \dotsc, u_i)$ of positive integers,
and any integers $e \geq 0$ and $t \geq 1$, we have: \begin{align*}
0 :_{F^e(M / (x_1^{u_1}, \dotsc, x_i^{u_i}) M)} x_{i+1}^t
&= 0 :_{F^e(M) / (x_1^{u_1 q}, \dotsc, x_i^{u_i q}) F^e(M)} x_{i+1}^t \\
&\subseteq 0^*_{F^e(M) / (x_1^{u_1 q}, \dotsc, x_i^{u_i q}) F^e(M)} \\
&= 0^*_{F^e(M / (x_1^{u_1}, \dotsc, x_i^{u_i}) M)}.
\end{align*}
The containment follows from setting $t_j = u_j q$ for $1 \leq j \leq i$ and $t_{i+1} = t$
in (\ref{eq:pdepth+def}) above.  To see that (b) $\Rightarrow$ (c),
for each $i$ and $e$, simply set $u_j = 1$ for $1 \leq j \leq i$ and $t = p^e$
in the definition of phantom $M$-regularity to obtain the defining containment for ghost
$M$-regularity.

To see that (c) $\Rightarrow$ (d), it suffices to show that if an element $x$ is ghost $M$-regular, 
then $x^q$ is $G^e(M)$-regular for all $e$.  So fix $e$ and let 
$z \in F^e(M)$ such that $x^q \bar{z} = \bar{0}$
in $G^e(M)$.  That is, $x^q z \in 0^*_{F^e(M)}$.  Then for all $q' \gg 0$, $c x^{q q'} z^{q'} = 0$,
so that by ghost $M$-regularity of $x$, $c z^{q'} \in 0^*_{F^{e+e'}(M)}$, whence $c_1 z^{q' q_0} = 0$
for all such $q'$, and hence $z \in 0^*_{F^e(M)}$, \emph{i.e.}, $\bar{z} = \bar{0}$.  So $x^q$
is $G^e(M)$-regular for all $e$.

The implication (d) $\Rightarrow$ (e) is trivial.

To show (e) $\Rightarrow$ (c), it suffices to show that if an element $x$ has the property that
$x^q$ is $G^e(M)$-regular for $e \gg 0$, then $x$ is ghost $M$-regular.  So fix $e$ and let $z \in F^e(M)$
such that $x^q z = 0$.  Then for all $q'$, $x^{q q'} z^{q'} = 0$, and for sufficiently large $q'$
the hypothesis then shows that $z^{q'} \in 0^*_{F^{e+e'}(M)}$, and therefore $c z^{q'q_0} = 0$ for all such
$q'$.  Hence, $z \in 0^*_{F^e(M)}$, which completes the proof that $x$ is ghost $M$-regular.

It remains to show that (c) $\Rightarrow$ (a).  Accordingly,
we assume $\bfx$ is a ghost $M$-regular sequence.

If $e \geq 0$, $n \geq 1$, and $t_1, \dotsc, t_n \geq 1$ are integers, 
then by the division algorithm, we have $t_1 = d q + r$, where $1 \leq r \leq q$.
We assume by induction on $n$ that for any $M$, any ghost $M$-regular sequence of length
$n-1$ satisfies (a), and then we prove by 
induction \emph{on $d$} that for any $e$ and any $n$-tuple $t_1, \dotsc, t_n$
of positive integers where $t_1 = d q + r$, where $1 \leq r \leq q$,
(\ref{eq:pdepth+def}) holds.

First take the case where $n=1$ and $d=0$.\footnote{This case is due to
a hint from Ian Aberbach.\cite{Abpers}}  That is, $x_1^{t_1} z = 0$, where
$t_1 \leq q$ and $z \in F^e(M)$.  Then $x_1^q z = x_1^{q-t_1} x_1^{t_1} z
= 0$, so that by ghost $M$-regularity of $x_1$, $z \in 0^*_{F^e(M)}$.
If $d>0$, then $t_1 > q$, and
$x_1^{t_1 - q} x_1^q z = x_1^{t_1} z = 0$, so that by induction on $d$,
$x_1^q z \in 0^*_{F^e(M)}$.  Then for any $e' \geq e_0$, $x_1^{q q'} (c z^{q'})
= c x_1^{q q'} z^{q'} = c (x_1^q z)^{q'} = 0$ in $F^{e+e'}(M)$.  Then by ghost $M$-regularity of 
$x_1$, $c z^{q'} \in 0^*_{F^{e+e'}(M)}$.  Hence, $c_1 z^{q'q_0} = 
c(c z^{q'})^{q_0} = 0$, and since this holds for all $e' \geq e_0$, it follows
that $z \in 0^*_{F^e(M)}$.

Now consider the case where $n>1$.  Suppose first that $d=0$, so that $t_1 \leq 
q$.  If $t_1 = q$, then we 
have \[
x_n^{t_n} z_n \in (x_1^q, x_2^{t_2}, \dotsc, x_{n-1}^{t_{n-1}}) F^e(M)
 = (x_1 M)^{[q]}_M + (x_2^{t_2}, \dotsc, x_{n-1}^{t_{n-1}}) F^e(M).
\]
Then, letting $\overline{z_n}$ be the image of
$z_n$ in $F^e(M / x_1 M) = F^e(M) / (x_1 M)^{[q]}_M = F^e(M) / (x_1^q F^e(M))$, we have \[
x_n^{t_n} \overline{z_n} \in (x_2^{t_2}, \dotsc, x_{n-1}^{t_{n-1}}) F^e(M / x_1 M),
\]
so that since $x_2, \dotsc, x_n$ satisfies (a) on the module $M / x_1 M$
(by induction on $n$), \begin{align*}
\overline{z_n} &\in \left( (x_2^{t_2}, \dotsc, x_{n-1}^{t_{n-1}}) F^e(M / x_1 M) \right)^*_{F^e(M/x_1 M)} \\
&= \frac{\left( x_1^q F^e(M) + (x_2^{t_2}, \dotsc, x_{n-1}^{t_{n-1}}) F^e(M) \right)^*_{F^e(M)}}{x_1^q F^e(M)} \\
&= \frac{\left( (x_1^q, x_2^{t_2}, \dotsc, x_{n-1}^{t_{n-1}}) F^e(M) \right)^*_{F^e(M)}}{x_1^q F^e(M)},
\end{align*}
which means that 
\[
z_n \in \left( (x_1^q, x_2^{t_2}, \dotsc, x_{n-1}^{t_{n-1}}) F^e(M) 
\right)^*_{F^e(M)}.
\]

Hence, we may assume that $t_1 < q$.  Say $x_1^{t_1} z_1 + \dotsc + x_n^{t_n} 
z_n = 0$.  Setting $v = q - t_1$ and multiplying the equation by $x_1^v$, we have \[
x_1^q z_1 + x_2^{t_2} (x_1^v z_2) + \dotsc + x_n^{t_n} (x_1^v z_n) = 0.
\]
By the case $t_1 = q$ above, $x_1^v z_n \in \left( (x_1^q, x_2^{t_2}, \dotsc, x_{n-1}^{t_{n-1}}) F^e(M) 
\right)^*_{F^e(M)}$, so that for any $e' \geq e_0$, there exist $y_1, \dotsc, y_{n-1}
\in F^{e+e'}(M)$ such that \[
c x_1^{v q'} z_n^{q'} = x_1^{q q'} y_1 + \sum_{i=2}^{n-1} x_i^{t_i q'} y_i.
\]
Collecting like terms, we have \begin{equation}\label{eq:col1}
x_1^{v q'} (x_1^{t_1 q'} y_1 - c z_n^{q'}) + \sum_{i=2}^{n-1} x_i^{t_i q'} y_i = 0.
\end{equation}
However, since ghost $M$-regular sequences are permutable, $x_2, \dotsc, x_{n-1}, x_1$
is a ghost $M$-regular sequence, so by induction on $n$, $x_2, \dotsc, x_{n-1}, x_1$
satisfies (a) for $M$.  Then Equation~\ref{eq:col1} gives that \[
x_1^{t_1 q'} y_1 - c z_n^{q'} \in \left( (x_2^{t_2 q'}, \dotsc, x_{n-1}^{t_{n-1}q'}) 
F^{e+e'}(M) \right)^*_{F^{e+e'}(M)}
\]
Taking to the $q_0$'th power and multiplying by $c$, we have \[
c x_1^{t_1 q' q_0} y_1^{q_0} - c_1 z_n^{q'q_0} \in \left((x_2^{t_2}, \dotsc, 
x_{n-1}^{t_{n-1}}) F^e(M)\right)^{[q'q_0]}_{F^e(M)}.
\]
Hence, $z_n \in \left((x_1^{t_1}, \dotsc, x_{n-1}^{t_{n-1}}) F^e(M)
\right)^*_{F^e(M)}.$

We are finally ready for the case where $n>1$ and $d>0$.  Say \begin{equation}\label{eq:6'}
x_1^{t_1} z_1 + \dotsc + x_n^{t_n} z_n = 0.
\end{equation}
We have $x_1^{t_1} z_1 = x_1^{t_1 - q} (x_1^q z_1)$, so that the
above equation yields, by induction on $d$, that
$z_n \in ((x_1^{t_1 - q}, x_2^{t_2}, \dotsc, x_{n-1}^{t_{n-1}})
F^e(M))^*_{F^e(M)}$.  Hence, for any $e' \geq e_0$, there exist
$y_1, \dotsc, y_{n-1} \in F^{e+e'}(M)$ such that \begin{equation}
\label{eq:y1final}
c z_n^{q'} = x_1^{(t_1 - q) q'} y_1 + \sum_{i=2}^{n-1} x_i^{t_i q'} y_i.
\end{equation}
Now, after applying the operator $c(-)^{[q']}$ to Equation~\ref{eq:6'} and
multiplying Equation~\ref{eq:y1final} by $x_n^{t_n q'}$, we can combine them
after solving each for the term $c x_n^{t_n q'} z_n^{q'}$, to get:
\begin{equation}\label{eq:9}
x_1^{(t_1 - q) q'} (c x_1^{q q'} z_1^{q'} + x_n^{t_n q'} y_1)
+ \sum_{i=2}^{n-1} x_i^{t_i q'} (c z_i^{q'} + x_n^{t_n q'} y_i) = 0.
\end{equation}
Since $x_2, \dotsc, x_{n-1}, x_1$ is a ghost
$M$-regular sequence, by permutability, hence one that satisfies (a) on $M$
by induction on $n$, Equation~\ref{eq:9} yields \[
c x_1^{q q'} z_1^{q'} + x_n^{t_n q'} y_1 \in 
\left( ((x_2^{t_2}, \dotsc, x_{n-1}^{t_{n-1}})F^e(M))^{[q']}_{F^e(M)}
\right)^*_{F^{e+e'}(M)}.
\]
Hence,  $c x_n^{t_n q' q''} y_1^{q''} \in \left( (x_1^qq', x_2^{t_2 q'}, \dotsc,
x_{n-1}^{t_{n-1} q'}) F^e(M) \right)^{[q'']}_{F^{e+e'}(M)}$ for any $e'' \geq e_0$.
But by induction on $n$, $x_2, \dotsc, x_n$ satisfies (a) for $(M / x_1 M)$, so we have \[
c y_1^{q''} \in \left( \left( (x_1^{q q'}, x_2^{t_2 q'}, \dotsc,
x_{n-1}^{t_{n-1} q'}) F^{e+e'}(M) \right)^{[q'']}_{F^{e+e'}(M)}
\right)^*_{F^{e+e'+e''}(M)}.
\]
Since this holds for any $e'' \geq e_0$ and $R$ has a weak test element, it
follows that $y_1 \in \left( (x_1^{q q'}, x_2^{t_2 q'}, \dotsc, x_{n-1}^{t_{n-1} q'})
F^{e+e'}(M) \right)^*_{F^{e+e'}(M)}.$  But then by Equation~\ref{eq:y1final},
\begin{multline*}
c z_n^{q'} \in x_1^{(t_1 - q) q'} \left( (x_1^{q q'}, x_2^{t_2 q'}, \dotsc, x_{n-1}^{t_{n-1} q'})
 F^{e+e'}(M) \right)^*_{F^{e+e'}(M)} \\
 + (x_2^{t_2 q'}, \dotsc, x_{n-1}^{t_{n-1} q'}) F^{e+e'}(M)
\end{multline*}
\[
\subseteq \left( \left( (x_1^{(t_1-q) + q}, x_2^{t_2}, \dotsc, x_{n-1}^{t_{n-1}})
F^{e}(M) \right)^{[q']}_{F^e(M)} \right)^*_{F^{e+e'}(M)}.
\]
Since this holds for all $e' \geq e_0$, and $R$ has a weak test element,
it follows that $z_n \in \left((x_1^{t_1}, \dotsc, x_{n-1}^{t_{n-1}}) F^e(M)
\right)^*_{F^e(M)}$.
\end{proof}

In Theorem~\ref{thm:flex} we need the following technical condition.
\begin{definition}
Let $R$ be a Noetherian ring of characteristic $p>0$.  We say that a
finitely generated $R$-module $M$ \emph{satisfies the avoidance
condition} (or it \emph{satisfies avoidance}) if for 
any quotient module $N$ of $M$ and any ideal $I 
\subseteq R$ such that \[
I \subseteq \bigcup \bigcup_{e \geq 0} \func{Ass} G^e(N),
\]
there is some $e \geq 0$ and some $\p \in \func{Ass} G^e(N)$ such that
$I \subseteq \p$.
\end{definition}

\begin{av}\label{av}
It may well be that every finitely generated module over any such
$R$ satisfies avoidance.  Certainly it is known when the ring itself
satisfies \emph{countable prime avoidance}\footnote{Countable prime avoidance
is the condition that any ideal 
which is contained in countable union of primes must be contained in one of the primes.}, 
which is the case if the ring is complete \cite[Lemma 3]{BurchCPA} or if it contains an
uncountable field \cite[Remark 2.17, for instance]{HHexponent}.  The set $\cup_{e \geq 0} 
\ffunc{Ass}_R  G^e_R(N)$ has a tendency to be finite, and if 
all quotients
of $M$ have the property that the corresponding union for the
quotient is finite, then the usual prime avoidance property
yields that $M$ satisfies avoidance.  However, Singh and Swanson
\cite{SiSwUAss} have given an example of a Noetherian normal 
ring $R$ where every 
ideal is tightly closed, along with an ideal in that ring where
the union of the primes associated to all the Frobenius powers (hence
also all the tight closures of the Frobenius powers) of 
the ideal is not a finite set.  Hence, we cannot rely on the
union to be actually finite.
\end{av}

The following is the main theorem of this section, and will be 
proved after Lemma~\ref{lem:cohoind}.  The corollary following its proof answers in the 
affirmative, at least for modules satisfying avoidance, the 
question raised in \cite{AbPPD}: All phantom $M$-regular sequences
in an ideal $I$ have the same length.

\begin{thm}\label{thm:flex}\footnote{See also Corollary~\ref{cor:flexrigid}.}
Let $R$ be a Noetherian local ring of characteristic $p>0$ which
contains a weak test element $c$, and let $M$ be a finitely generated
$R$-module.  Take any sequence $\mathbf{x} = x_1, \dotsc, x_n$ of 
elements of $\m$ and any positive integer $d \leq n$.  Consider
the following conditions \begin{enumerate}
\item[(a)] The ideal $(\mathbf{x})$ contains a
phantom $M$-regular sequence of length $d$.
\item[(b)] $K.(\bfx; M)$ has stably phantom homology at $n-i$ for
$i = 0, 1, \dots, d-1$.
\end{enumerate}
Then (a) $\Rightarrow$ (b), and if $M$ satisfies avoidance then
(b) $\Rightarrow$ (a) as well.
\end{thm}

We prove this theorem through a series of lemmas, the first of which
provides a characterization of phantom zerodivisors:

\begin{lemma}\label{lem:phass}\cite[essentially Lemma 1.1]{nmepdepBC}
Let $R$ be a Noetherian local ring of characteristic $p>0$ and $M$ a 
finitely generated $R$-module.  Then the set of phantom zerodivisors for 
$M$ is contained in $\bigcup
\bigcup_e \func{Ass} G^{e}(M)$.  If $R$ has a $q_0$-weak test element $c$, then every
element of this union is a phantom zerodivisor for $M$.
\end{lemma}

\begin{proof}
For the first part, let $x$ be a phantom zerodivisor for $M$.  Then
there is some $e \geq 0$ and $t > 0$ such that $0 :_{F^{e}(M)} x^{t} \not
\subseteq 0^{*}_{F^{e}(M)}$.  That is, there is some $z \in F^{e}(M)
\setminus 0^{*}_{F^{e}(M)}$ with $x^{t}z = 0$.  Then $x^{t}
\overline{z} = \overline{0}$ in $G^{e}(M)$, where $\overline{z} \neq \bar{0}$,
so there is some $\p \in
\Ass G^{e}(M)$ with $x^{t} \in \p$.  Since $\p$ is prime and thus 
radical, $x \in \p$.

For the second part, let $x \in \p$ for some $\p \in \Ass G^{e}(M)$
for some $e$.  Then there is some $z \in F^{e}(M)$, $z \not \in
0^{*}_{F^{e}(M)}$, with $\p = \overline{0} :_{G^{e}(M)} \overline{z}$,
which means that $x z \in 0^{*}_{F^{e}(M)}$.  Then for all large powers
$q' \gg 0$ of $p$, $c x^{q'}z^{q'} = 0$.  If $x$ is phantom $M$-regular,
this together with the last equation implies that $c z^{q'} 
\in 0^*_{F^{e+e'}(M)}$.  Hence, $c_1 z^{q'q_0} = 
c (c z^{q'})^{q_0} = 0$, from which we conclude that $z \in 0^*_{F^e(M)}$,
contrary to assumption.  Thus, $x$ is a phantom zerodivisor for $M$.
\end{proof}

Next, we prove the case $d=1$ of Theorem~\ref{thm:flex}:

\begin{lemma}\label{lem:case1}
Let $R$ be a Noetherian local ring of characteristic $p>0$ containing
a $q_0$-weak test element $c$.  Let $M$ be a
finitely generated $R$-module and $\mathbf{x}$ a sequence in $\m$ of
length $n$.  Consider the following conditions \begin{enumerate}
\item[(a)] The ideal $(\mathbf{x})$ contains a phantom $M$-regular element.
\item[(b)] $K_.(\mathbf{x}; M)$ has stably phantom homology at $n$.
\end{enumerate}
Then (a) $\Rightarrow$ (b), and if $M$ satisfies
avoidance then (b) $\Rightarrow$ (a) as well.
\end{lemma}

\begin{proof}
Note that an element of $H_n(\mathbf{x}^{[q]}; F^{e}(M))$ is
precisely an element of $F^{e}(M)$ which is killed by each
$x_{i}^{q}$; $i = 1, \ldots, n$, and we have that $B_n(\bx{q}; F^e(M)) = 0$.
By Proposition~\ref{prop:killphantom}(b), we need to show that the existence of a phantom
$M$-regular element in $(\mathbf{x})$ implies (and if $M$ satisfies avoidance,
is equivalent to) the assertion
that for any $e\geq 0$, any element of $F^{e}(M)$ which is killed by all of
the $x_{i}^{q}$ is an element of $0^*_{F^{e}(M)}$.

First suppose that we have a phantom $M$-regular element $y \in
(\mathbf{x})$.  Take any nonnegative integer $e$, and let $z \in
F^{e}(M)$ be an element annihilated by all of the $x_{i}^{q}$.  Then
since $y^{q}$ is a linear combination of the $x_{i}^{q}$, $y^{q}z=0$,
so by definition of phantom (or ghost) $M$-regularity, $z \in 0^{*}_{F^{e}(M)}$.
Hence $K.(\mathbf{x};M)$ has stably phantom homology at $n$.

For the other direction, assume that $M$ satisfies avoidance and that
$(\mathbf{x})$ has no phantom
$M$-regular elements.  By Lemma~\ref{lem:phass} and avoidance, there is some $e$ 
and some $z \in F^e(M) \setminus 0^*_{F^e(M)}$ such that $\p = 0^*_{F^e(M)}
:_{F^e(M)} z$ is a prime ideal containing $(\mathbf{x})$.  Then for each $i = 1, \ldots, n$,
$x_{i}z \in 0^{*}_{F^e(M)}$, so $x_i^q z \in 0^*_{F^e(M)}$, and thus 
$x_i^{qq'}(c z^{q'}) = c(x_{i}^{q}z)^{q'}=0$, for all $e'\geq e_0$ and all $i$.
If $K.(\bfx; M)$ had stably phantom homology at $n$,
then since $c z^{q'}$ is killed by all the $x_i^{qq'}$, we would have 
$c z^{q'} \in 0^{*}_{F^{e+e'}(M)}$, hence $c_1 z^{q'q_0} =
c (c z^{q'})^{q_0} = 0$, for all $e'$.  Thus, $z \in 0^{*}_{F^{e}(M)}$,
which is a contradiction.  Hence, $K.(\bfx; M)$ does not have stably phantom homology at $n$.
\end{proof}

The final preparatory lemma provides the inductive step in Theorem~\ref{thm:flex}.  In preparation,
note the easy fact that 
\[
F^{e}(K.(\bfx; M)) \cong K.(\bx{p^e}; F^{e}(M))
\]
as complexes of $R$-modules, for any $e$.

\begin{lemma}\label{lem:cohoind}
Let $R$ be a Noetherian local ring of characteristic $p>0$ containing a $q_0$-weak test
element $c$, and let $M$ be a finitely generated $R$-module.  Let $y_1, \ldots, y_u$ be 
a phantom $M$-regular sequence in $(\mathbf{x})$, where
$u \geq 1$ and $\mathbf{x} = x_1, \ldots, x_n$.  Then
the following are equivalent: \begin{enumerate}
\item[(a)] $K.(\mathbf{x}; M)$ has stably phantom homology at $n-i$ for $i = 0, 1, \ldots, u$.
\item[(b)] $K.(\mathbf{x}; M  / y_1 M)$ has stably phantom homology at $n-j$ for $j=0, 1, \ldots, u-1$.
\item[(c)] $K.(\mathbf{x}; M / (y_1, \ldots, y_u) M)$ has stably phantom homology at $n$.
\end{enumerate}
\end{lemma}

\begin{proof}
\noindent
For parts (a) and (b), note that the following sequence of complexes: \[
0 \rightarrow K.(\mathbf{x}; M) \stackrel{\alpha.}{\rightarrow} K.(\mathbf{x}; M)
\stackrel{\beta.}{\rightarrow} K.(\mathbf{x}; M / y_1 M) \rightarrow 0
\]
satisfies the conditions of Proposition~\ref{prop:pseudolesABBR}, where for each
$i$, $\alpha_i$ is multiplication by $y_1$ and $\beta_i$ is the canonical surjection
associated to $\alpha_i$, since $\ker (y_1^q: F^e(M) \rightarrow F^e(M)) \subseteq 0^*_{F^e(M)}$
for all $e \geq 0$.  Also note that for this $\alpha.$, we have $H_{i}(F^e(\alpha.)) = 0$
for all $i$ and all $e \geq 0$, since $y_1^q \in (\bx{q})$ for all $q$. (see \cite[Proposition 1.6.5]{BH})

Suppose that (a) holds.  Then for $j = 0, \ldots, u-1$, 
$K.(\mathbf{x}; M)$ has stably phantom homology at $n-j$ and 
$n-(j+1) = (n-j)-1$, so that by Proposition~\ref{prop:pseudolesABBR}(\ref{it:LMph}), 
$K.(\bfx; M / y_1 M)$ has stably phantom homology at $n-j$.

Conversely, suppose that (b) holds.  First take $i \in \{ 0, \ldots, u-1 \}$,
and $[z] \in H_{n-i}(\bx{q}; F^e(M))$.  Then by 
Proposition~\ref{prop:pseudolesABBR}(\ref{it:Nph}b), $[c z^{q'}] \in \im H_{n-i}(F^{e+e'}(\alpha.))
= 0$ for all $e' \geq e_0$, since $K.(\bfx; M / y_1 M)$ has stably phantom homology at $n-i$.
Hence $[z]$ is a phantom element of $H_{n-i}(\bx{q}; F^e(M))$.
On the other hand, take $[w] \in H_{n-u}(\bx{q}; F^e(M))$.  Then since $[\alpha_{n-u}(w)] = 0$ and 
$K.(\bfx; M / y_1 M)$ has stably phantom homology at $n-u+1$, it follows from
Proposition~\ref{prop:pseudolesABBR}(\ref{it:Nph}a), that $[w]$ is a phantom element of $H_{n-u}(\bx{q};
F^e(M))$.  Hence, (a) holds.

Finally, we prove that $(b) \Leftrightarrow (c)$:
If $u=1$, then the two expressions are identical.  If $u >1$, then since 
$y_2, \ldots, y_u$ is a phantom $(M/y_1 M)$-regular sequence, the equivalence
follows by induction on $u$ and the fact that parts (a) and (b) are equivalent.
\end{proof}

\begin{proof}[Proof of Theorem~\ref{thm:flex}.]
We proceed by induction on $d$.  The case $d=1$ is
Lemma~\ref{lem:case1}, so we may assume that $d \geq 2$ and that we've
proved the theorem for smaller $d$.

First suppose that $(\mathbf{x})$ contains a phantom $M$-regular
sequence $y_1, \ldots, y_d$.  Then
$y_d$ is phantom $\left(M / (y_1, \ldots, y_{d-1}) M\right)$-regular, 
so by Lemma~\ref{lem:case1}, $K.(\bfx; M / (y_1, \ldots, y_{d-1})M)$
has stably phantom homology at $n$, whence by Lemma~\ref{lem:cohoind}, 
$K.(\bfx; M)$ has stably phantom homology at $n-i$ for all $i = 0, 1, 
\ldots, d-1$.

For the converse, suppose that $K.(\bfx; M)$ has stably phantom homology at $n-i$ for
$i = 0, 1, \ldots, d-1$, and that $M$ satisfies
avoidance.  By induction on
$d$, there is a phantom $M$-regular sequence $y_1, \ldots, y_{d-1}$
in $(\mathbf{x})$.  Therefore, by Lemma~\ref{lem:cohoind},
$K.(\bfx; M / (y_1, \ldots, y_{d-1})M)$ has stably phantom homology at $n$,
and then we use Lemma~\ref{lem:case1} to guarantee the existence of
a phantom $(M / (y_1, \ldots, y_{d-1})M)$-regular element $y_d$, giving the
desired phantom $M$-regular sequence $y_1, \ldots, y_d$ in $(\mathbf{x})$.
\end{proof}

\begin{cor}\label{cor:flex}
Let $R$ be a Noetherian local ring of characteristic $p>0$ which
contains a weak test element $c$, let $M$ be a finitely generated
$R$-module with avoidance, and let $I$ be a proper
ideal of $R$.  Then all maximal phantom $M$-regular sequences in $I$ have
the same length: $\pdp{I} M$.
\end{cor}
\begin{proof}
Let $\mathbf{y} = y_1, \dotsc, y_t$, be a maximal phantom $M$-regular sequence in $I$, and
let $d = \pdp{I} M$.  Let $\mathbf{x} = x_1, \dotsc, x_n$ be a generating set for $I$.  By 
Lemma~\ref{lem:case1}, since $I$ does not contain any phantom 
$(M / \mathbf{y}M)$-regular elements, 
$K.(\bfx; M / \mathbf{y}M)$ does not have stably phantom homology at $n$.  Then by 
Lemma~\ref{lem:cohoind}, there is some $i \in \{0, \dotsc, t\}$
such that $K.(\bfx; M)$ does not have stably phantom homology at $n-i$.
On the other hand, by Theorem~\ref{thm:flex}, since $I$ contains a phantom 
$M$-regular sequence of length $d$, 
$K.(\bfx; M)$ has stably phantom homology at $n-j$ for $j = 0, \dotsc, d-1$.
Hence, $t \geq i \geq d$, but we also know that $d \geq t$ by definition
of phantom depth, so
$d=t$.  Thus, all maximal phantom $M$-regular sequences in $I$ have length $d$.
\end{proof}

\section{Stable phantom rigidity of Koszul complexes and ghost \emph{M}-regular sequences}\label{sec:phseqK}
In this section we show (Theorem~\ref{thm:stphrigidK}) that over a prime characteristic 
Noetherian local ring with a weak test element,
any Koszul complex is ``stably phantom rigid.''
According to Auslander \cite{Ausunram}, a complex $K.$ is called \emph{rigid} if whenever 
$L$ is finitely generated and $i>0$ an integer such that
$\Tor_i(K.,L) = 0$, it follows that $\Tor_j(K.,L) = 0$ for
all $j>i$.  A module is called \emph{rigid} if its minimal free resolution
is a rigid complex.  

In \cite{ABcodmul}, Auslander and Buchsbaum showed that the Koszul complex on any finite
set of elements of a ring is rigid.
Auslander then used this along with Serre's diagonalization argument to show \cite[Theorem 2.1]{Ausunram} 
that any finitely generated module over an equicharacteristic or unramified regular local ring is rigid.  Lichtenbaum
\cite{LicTor} showed the same result for ramified regular local rings.  If we pass to non-regular
local rings, examples of non-rigid modules of infinite projective dimension are easily constructed.
However, Auslander's question in \cite{Ausunram}, of whether modules of finite projective dimension
over a Noetherian local ring are always rigid, was open for over 30 years.
As Hochster noted in \cite[Chapter 2]{HoCBMS}, if it were true, it would yield several
interesting results which were at that time conjectures.  Since then, these consequences
of rigidity were shown to be true, but Heitmann showed Auslander's ``rigidity conjecture'' is false.\cite{Heitrigid}
The ring in Heitmann's counterexample is not a complete intersection, and indeed the rigidity conjecture
is open for complete intersections.  For work on rigidity over complete intersections,
see \cite{MuRLR}, \cite{HWrigid}, \cite{HWrigid2}, \cite{Jocx}, and \cite{HJW}.

In this section, we provide an analogue in Theorem~\ref{thm:stphrigidK}, in terms of stable phantomness of homology,
of Auslander and Buchsbaum's result on rigidity of the Koszul complex.  

The same method is used to obtain 
a Koszul homology criterion for a sequence of 
elements to be a ghost $M$-regular sequence (Theorem~\ref{thm:ghseqK}).  Aberbach also gives a Koszul
homology criterion for a sequence to be phantom $M$-regular
in \cite[Theorem 3.3.8]{AbPPD}, so by Proposition~\ref{prop:gh=ph}, we already have
a Koszul homology criterion for a sequence to be ghost $M$-regular.  However, the
criterion given here seems easier to check, and in any case the stable phantom rigidity
result is interesting in its own right.  First note the following

\begin{lemma}\label{lem:xx'}
Let $(R,\m)$ be a Noetherian local ring of characteristic $p>0$ containing
a $q_0$-weak test element $c$. Let $\bfx = x_1, \dotsc, x_n$ be a sequence of
elements of $\m$, where $n \geq 1$, 
and let $\mathbf{x}' = x_1, \dotsc, x_{n-1}$.  Then for any $i>0$, if 
$K.(\bfx; M)$ has stably 
phantom homology at $i$, so does $K.(\bfx'; M)$.
\end{lemma}
\begin{proof}
Let $[y] \in H_i(\bxp{q}; F^e(M))$ for some fixed $e \geq 0$.  Recall (e.g.
\cite[Corollary 1.6.13(a)]{BH}) that we
have an exact sequence: \begin{equation}\label{eq:xnexact}
H_i(\bxp{q}; F^e(M)) \arrow{\pm x_n^q} H_i(\bxp{q}; F^e(M)) \arrow{g_e}
H_i(\bx{q}; F^e(M)).
\end{equation}
Then $g_e([y]) \in H_i(\bx{q}; F^e(M))$ is, of course, phantom, 
so that $g_{e+e'}([cy^{q'}])$ is $[0]$ for any $e' \geq e_0$.  Fix some
such $e'$.  Then by exactness of (\ref{eq:xnexact}) for $e+e'$, $[c y^{q'}] = 
x_n^{q q'} [w]$ for some $[w] \in H_i(\bxp{q q'}; F^{e+e'}(M))$.  That
is, there exist $v \in K_{i+1}(\bxp{q q'}; F^{e+e'}(M))$ and $w
\in K_i(\bxp{q q'}; F^{e+e'}(M))$ such that $\partial(w) = 0$ and \[
c y^{q'} = \partial(v) + x_n^{q q'} w.
\]
Since $y \in Z_i(\bxp{q}; F^e(M))$ was arbitrarily chosen and since $x \in \m$, 
we have that \begin{multline}\label{eq:8'}
c Z_i(\bxp{q}; F^e(M))^{[q']}_{K_i(\bxp{q}; F^e(M))} \\
 \subseteq B_i(\bxp{q}; F^e(M))^{[q']}_{K_i(\bxp{q}; F^e(M))} \\
 + \m^{[q q']} Z_i(\bxp{q q'}; F^{e+e'}(M)).
\end{multline}
Since (\ref{eq:8'}) holds for all $e$ and all $e' \geq 0$, Proposition~\ref{prop:Nak*a}
implies that $Z_i(\bxp{q}; F^e(M)) \subseteq B_i(\bxp{q}; F^e(M))^*_{K_i(\bxp{q}; F^e(M))}$
for all $e \geq 0$.  That is, $K.(\bfx; M)$ has stably phantom homology at $i$.
\end{proof}

\begin{thm}[Stable phantom rigidity of Koszul complexes]\label{thm:stphrigidK}
Let $(R,\m)$ be a Noetherian local ring with a $q_0$-weak test element $c$, and
let $M$ be a finitely generated $R$-module.  Let $\bfx = x_1, \dotsc, x_n$ be
any sequence of elements of $\m$, and let $i \geq 1$ be an integer.  
If $K.(\bfx; M)$ has stably phantom homology at $i$, then it has stably phantom homology at
$j$ for all $j \geq i$.
\end{thm}

\begin{proof}[Proof of theorem]
It suffices to prove the result for $j = i+1$, which we shall do by induction on $n$.  If 
$n=1$, then there is nothing to prove because $H_j(x_1^q; F^e(M)) = 0$ whenever
$j >1$.

So assume that $n>1$ and that the result has been shown for sequences of length $n-1$.
Note that for any $i \geq 1$, the following equation for differentials of Koszul
complexes holds, if we set $\bfxp$ to be the sequence $x_1, \dotsc, x_{n-1}$ and
$\delta_{j,\mathbf{y}}: K_j(\mathbf{y}; M) \rightarrow K_{j-1}(\mathbf{y}; M)$
is the $j$'th differential in the Koszul complex:
\begin{equation}\label{eq:8''}
\delta_{i+1, \bfx}^{[q]} = \left(\begin{matrix}
 \delta_{i+1, \bfxp}^{[q]} & (-1)^{n-i+1} x_n^q \\
 0 & \delta_{i,\bfxp}^{[q]}
\end{matrix}\right),
\end{equation}
using the natural identification \[
K_{i+1}(\bx{q}; F^e(M)) \simeq K_{i+1}(\bxp{q}; F^e(M))\oplus K_i(\bxp{q}; F^e(M))
\]
Let $\left(\begin{matrix}
u \\
v \end{matrix}\right) \in \ker \delta_{i+1,\bfx}^{[q]}$.  Then by (\ref{eq:8''}),
we have:
\begin{equation}\label{eq:uv}
\delta_{i+1,\bfxp}^{[q]}(u) = (-1)^{n-i} x_n^q v
\end{equation}
and \begin{equation}\label{eq:9'}
\delta_{i,\bfxp}^{[q]}(v) = 0.
\end{equation}
$K.(\bfxp; M)$ has stably phantom homology at $i$ by Lemma~\ref{lem:xx'}, which combines with (\ref{eq:9'})
to show that $v \in (\im \delta_{i+1, \bfxp}^{[q]})^*_{K_i(\bxp{q}; F^e(M))}.$
Fixing $e' \geq e_0$, there exists $y \in K_{i+1}(\bxp{q q'}; F^{e+e'}(M))$
such that \begin{equation}\label{eq:vy}
c v^{q'} = \delta_{i+1, \bfxp}^{[q q']}(y).
\end{equation}

Combining Equations \ref{eq:uv} and \ref{eq:vy}, we have \[
\delta_{i+1, \bfxp}^{[q q']}(c u^{q'}) = (-1)^{n-i} x_n^{q q'} c v^{q'}
= \delta_{i+1, \bfxp}^{[q q']}( (-1)^{n-i} x_n^{q q'} y).
\]
Thus, \[
c u^{q'} + (-1)^{n-i+1} x_n^{q q'} y \in \ker \delta_{i+1, \bfxp}^{[q q']}
\subseteq \left( \im \delta_{i+2, \bfxp}^{[q q']} \right)^*,
\]
where the last containment is by inductive hypothesis, since $\bfxp$ is a
sequence of length $n-1$.  Then there exists $z$
such that $c_1 u^{q' q_0} + (-1)^{n-i+1} x_n^{q q' q_0} c y^{q_0} = 
\delta_{i+2, \bfxp}^{[q q' q_0]}(z)$.  Combining everything together, 
we have that \begin{align*}
c_1 \left(\begin{matrix} u \\ v \end{matrix} \right)^{q' q_0}
&= \left(\begin{matrix}
\delta_{i+2, \bfxp}^{[q q' q_0]} & (-1)^{n-i} x_n^{q q' q_0} \\
0 & \delta_{i+1, \bfxp}^{[q q' q_0]}
\end{matrix}\right) \left(\begin{matrix} z \\ c y^{q_0} \end{matrix} \right) 
= \delta_{i+2, \bfx}^{[q q' q_0]} \left( \begin{matrix}
z\\ c y^{q_0}
\end{matrix} \right) \\
&\in \left( \im \delta_{i+2, \bfx}^{[q]} \right)^{ [q' q_0]}_{K_{i+1}
(\bx{q}; F^e(M))}.
\end{align*}
Thus, $\left(\begin{matrix} u \\ v \end{matrix}\right) \in (\im 
\delta_{i+2,\bfx}^{[q]} )^*_{K_{i+1}(\bx{q}; F^e(M))}$, showing that
$K_.(\bfx; M)$ has stably phantom homology at $i+1$.
\end{proof}

\begin{cor}\label{cor:flexrigid}
In Theorem~\ref{thm:flex}, one can replace the phrase ``at $n-i$ for
$i = 0, \dotsc, d-1$'' by the phrase ``at $n-d+1$''.  That is, we need only check
for stably phantom homology at one spot.
\end{cor}

\begin{thm}[Koszul homology characterization of ghost $M$-regular sequences]\label{thm:ghseqK}
Let $(R,\m)$ be a Noetherian local ring with a weak test element $c$, and
let $M$ be a finitely generated $R$-module.  Let
$\bfx = x_1, \dotsc, x_n$ be any sequence of elements of $\m$.  Then the
following conditions are equivalent: \begin{enumerate}
\item[(a)] $\bfx = x_1, \dotsc, x_n$ is a ghost $M$-regular sequence,
\item[(b)] $K.(\bfx; M)$ has stably phantom homology at 1,
\item[(c)] $K.(\bfx; M)$ has stably phantom homology at $j$ for all $j \geq 1$.
\end{enumerate}
\end{thm}

\begin{proof}
If $\bfx$ is a ghost $M$-regular sequence, then certainly the ideal $(\bfx)$
contains a ghost $M$-regular sequence of length $n$, from which Theorem~\ref{thm:flex}
gives us that (a) implies (c).  (c) implies 
(b) trivially.

Suppose that (b) is true.  We will show (a) by induction on $n$.
If $n=1$, then the result follows directly
from the definition of ghost $M$-regularity.  So assume $n>1$ and that the 
theorem holds for sequences of length $n-1$.  By Lemma~\ref{lem:xx'}, $K.(\bfxp; M)$ has stably
phantom homology at 1, where $\bfxp = x_1, \dotsc, x_{n-1}$, 
$\bfxp$ is a ghost $M$-regular
sequence by inductive hypothesis.  Let $z \in F^e(M)$ with $x_n^q z \in \bxp{q}F^e(M)$.  Let
$\delta_{r, \bfx}^{[q]}$
denote the differentials in the Koszul complex, as in the
proof of Theorem~\ref{thm:stphrigidK}.  Then 
there exist $z_1, \dotsc, z_n \in F^e(M)$ with $z = z_n$ and \[
\left(\begin{matrix}
z_1\\
z_2\\
\vdots\\
z_n
\end{matrix}\right) \in \ker \delta_{1, \bfx}^{[q]}
\subseteq (\im \delta_{2,\bfx}^{[q]})^*_{K_1(\bx{q}; F^e(M))},
\]
so that for any $e' \geq e_0$, \[
\left(\begin{matrix}
c z_1^{q'} \\
\vdots \\
c z_n^{q'} \end{matrix}\right) \in \im \delta_{2, \bfx}^{[q q']}.
\]
As in the proof of Theorem~\ref{thm:stphrigidK}, \[
\delta_{2,\bfx}^{[q q']} = \left( \begin{matrix}
\delta_{2,\mathbf{x}'}^{[q q']} & -x_n^{q q'} \\
0 & \delta_{1,\bfxp}^{[q q']}
\end{matrix}
\right),
\]
so \[
c z^{q'} = c z_n^{q'} \in \im \left( \begin{matrix}
0 & \delta_{1, \bfxp}^{[q q']}
\end{matrix} \right) = \left( (x_1^q, \dotsc, x_{n-1}^q) 
F^e(M) \right)^{[q']}_{F^e(M)},
\]
so that $z \in \left( (x_1^q, \dotsc, x_{n-1}^q) F^e(M) \right)^*_{F^e(M)}$,
whence $x_n$ is ghost $(M/ \bfxp M)$-regular, and so $\bfx$ is a ghost
$M$-regular sequence.
\end{proof}

In Aberbach's Koszul homology criterion, one simply replaces the $q$'th powers in 
Theorem~\ref{thm:ghseqK} with arbitrary powers.

\section{Comparison with minheight}
In \cite{HHphantom}, Hochster and Huneke define (and develop extensively) the notion of the minheight
of a module on an ideal.  By definition, if $I$ is an ideal of $R$, $M$ is a finitely generated
$R$-module, and $\p_1, \dotsc, \p_t$ are
the minimal primes of $M$, then the \emph{minheight of $M$ on $I$} (denoted $\ffunc{mnht}_I M$)
is defined to be $\min \{\func{ht} (I + \p_j) / \p_j \mid 1 \leq j \leq t \}$.  If $M = R$
we call it the \emph{minheight of $I$} and denote it by $\func{mnht} I$.
The following Proposition connects minheight with phantom depth.  

\begin{prop}\label{prop:minheight}
Let $(R,\m)$ be a Noetherian local ring of characteristic $p>0$, and
$J$ any proper ideal of $R$.
\begin{enumerate}
\item[(a)] If $R$ satisfies colon capturing (which holds, e.g., if $R$ is either a quotient of a 
Cohen-Macaulay local ring \cite[Theorem 3.1]{HuTC} or a module-finite extension of
a regular local ring \cite[Theorem 3.1A]{HuTC}), then $\func{mnht} J \leq \pdp{J} R$.
\item[(b)] If $R$ has a locally stable weak test element $c$ and $M$ is a finitely generated
$R$-module such that for every 
$\p \in \Spec R$, $M_\p$ is faithful as an $R_\p$-module, then $\ffunc{mnht}_J M \geq \pdp{J} M$.
\end{enumerate}
\end{prop}

The hypotheses of (b) hold when $R$ has a locally stable
weak test element and we take $M=R$, so in this case, $\func{mnht} J \geq \pdp{J} R$.  If, moreover, $R$ satisfies
colon capturing, then $\func{mnht} J = \pdp{J} R$, as 
Aberbach notes in the proof of \cite[Theorem 3.2.7]{AbPPD}.

We need to prove the following Lemma:

\begin{lemma}\label{lem:min}
Let $(R,\m)$ be a Noetherian local ring of characteristic $p>0$ with a
$q_0$-weak test element $c$, and let 
$M$ be a finitely generated faithful $R$-module.  Suppose that 
$y_1, \ldots, y_d$ is a phantom $M$-regular sequence in $\m$,
and $\q$ is a minimal prime of $M$.  Then $\m$ is not minimal over $\q +
(y_1, \ldots, y_{d-1})$.
\end{lemma}
\begin{proof}
Since $M$ is a faithful $R$-module, $\q$ is a minimal prime of $R$.  Then
$\q R_{\q}$ is a nilpotent ideal in $R_{\q}$, so there exists some $a \notin \q$ 
and some positive integer $t$ such that for all $q \geq t$, $a \q^q = 0$.

Assume that $\m$ is minimal over $\q + (y_1, \ldots, y_{d-1})$.
Then for some power $q'$ of $p$, 
$\m^{q'} \subseteq \q + (y_1, \ldots, y_{d-1})$.  Hence, for any $e$ so large that
$q = p^e \geq t$, \[
a y_d^{q q'} \in a \m^{[q q']} \subseteq (y_1^q, \ldots, y_{d-1}^q).
\]
Then for any $z \in F^e(M)$, \[
a y_d^{q q'} z \in (y_1, \dotsc, y_{d-1})^{[q]} F^e(M) = ((y_1, \dotsc, y_{d-1})M)^{[q]}_M.
\]

But since $y_1, \ldots, y_d$ is a phantom $M$-regular sequence, it follows
that \[
a z \in \left((y_1, \ldots, y_{d-1})M)^{[q]}\right)^*_{F^e(M)}.
\]
Hence, by the weak test element hypothesis, \[
c (a z)^{q_0} = c a^{q_0} z^{q_0} \in ((y_1, \ldots, y_{d-1}) M)^{[q q_0]}_M,
\]
so that since $z$ was an arbitrary element of $F^e(M)$ and since $F^{e+e_0}(M)$ is
generated as an $R$-module by all the $q_0$'th ``powers'' of such elements, it follows that \[
c a^{q_0} \in ((y_1, \ldots, y_{d-1})M)^{[q q_0]}_M :_R F^{e+e_0}(M)
= 0 :_R F^{e+e_0}(M / (y_1, \ldots, y_{d-1})M).
\]
Since $e$ can be taken to be arbitrarily large, and since $0 :_R F^e(N) \subseteq \m^{[q]}$
for any finitely generated $R$-module $N$, and any
$e \geq 0$,\footnote{Here is a proof of that fact.  Let \[
R^t \stackrel{\phi}{\rightarrow} R^s \rightarrow N \rightarrow 0
\]
be a minimal free presentation of $N$.  Then if we fix bases of $R^t$ and $R^s$,
all the entries of the matrix representing $\phi$ are in $\m$.  Moreover, 
$\phi^{[q]}$ is a minimal free presentation of $F^e(N)$, and all its entries
lie in $\m^{[q]}$.  Now take any $a \in 0 :_R F^e(N)$.  Then by exactness of
the presentation, for the first basis element (or any other basis element) $e_1$
of $R^s$, \[
a e_1^{[q]} \in \im \phi^{[q]} \subseteq \m^{[q]} R^s.
\]
Hence, $a \in \m^{[q]}$.} it now follows from the Krull
intersection theorem that $c a^{q_0} = 0 \in \q$.  However, $c \notin \q$ and
$a \notin \q$, which contradicts of the fact that $\q$ is prime.
\end{proof}

\begin{proof}[Proof of Proposition \ref{prop:minheight}]
For part (a), let $n = \func{mnht} J$.  Then $J$ contains a sequence $x_1, \ldots, x_n$ of
elements of $R$ such that for all $i = 1, \ldots, n$, $\func{mnht} (x_1,
\ldots, x_i) = i$.  Fix some integer $i$ with $0 \leq i < n$, and some 
sequence $t_1, \ldots, t_{i+1} \in \N_+$.  Suppose $y \in R$ such that
$x_{i+1}^{t_{i+1}} y \in (x_1^{t_1}, \ldots, x_i^{t_i})$.  Then by colon
capturing, $y \in (x_1^{t_1},
\ldots, x_i^{t_i})^*$.  Hence, $x_1, \ldots, x_n$ is a phantom $R$-regular
sequence, so $\ffunc{ph. depth}_J R \geq n = \func{mnht} J$.

For part (b), let $\q$ be a minimal prime of $M$, let 
$y_1, \ldots, y_d$ be a phantom $R$-regular sequence, where $d = \ffunc{ph. depth}_J
M$, and for some $0 \leq i<n$
let $\p$ be a minimal prime of $\q + (y_1, \ldots, y_i)$.  Then if $y_{i+1} \in \p$,
we have that $y_1/1, \ldots, y_{i+1}/1$ is a phantom $M_{\p}$-sequence, that $\q R_{\p}$
is a minimal prime of $M_{\p}$, that $M_\p$ is a faithful $R_\p$-module,
and that $\p R_{\p}$ is a minimal prime of $\q R_{\p} + (y_1/1, \ldots, y_i/1)$.  This
set of circumstances contradicts Lemma~\ref{lem:min}, so we conclude
that $y_{i+1} \notin \p$.  It follows from the criterion for minheight
in \cite[Proposition 2.2(c)]{HHphantom} that $\ffunc{ph. depth}_J M = d \leq \ffunc{mnht}_J M$.
\end{proof}

Thus, when $R$ has a weak test element and $M_\p$ is a faithful $R_\p$-module for all $\p$,
we have the following string of inequalities: \[
\ffunc{depth}_J M \leq \pdp{J} M \leq \ffunc{mnht}_J M \leq \ffunc{ht}_J M,
\]
and if $M = R$, the middle inequality is always an equality.
In particular, if $J = \m$, we have: \[
\func{depth} M \leq \func{ph. depth} M \leq \ffunc{mnht}_\m M \leq \dim M.
\]

\section{Phantom assassinators, stable primes, and the associated primes of reduced Frobenius powers}
\label{sec:phass}
There are three analogues of associated primes of a finitely generated $R$-module $M$ that
I know of in tight closure theory.  The first is the set \[
\GAss{M},
\]
which we have already seen at work many times in this paper.  The second is the \emph{phantom assassinator}
$\Assph{M}$.  The latter notion is explored in \cite{AHH}, but there it is only defined for
modules of finite phantom projective dimension, and is defined to consist of all primes $\p$ such
that $\ffunc{ppd}_{R_\p} M_\p = \func{mnht} \p R_\p$.  However, I see no reason to restrict this
notion to the case of finite ppd, and so I propose the following definition, which 
is equivalent to the one given in \cite{AHH} under mild conditions on $R$\footnote{In particular, 
they are equivalent whenever
Aberbach's phantom Auslander-Buchsbaum Theorem \cite[Theorem 3.2.7]{AbPPD} holds, which in turn is true
over nearly all rings of interest: see \cite[Proposition 5.4 and Discussion 5.5]{AHH}.  There
such rings are said to have \emph{acceptable type}.} when $\ffunc{ppd}_R M < \infty$: \[
\Assph M := \{ \p \in \Spec R \mid \pdp{\p R_\p} M_\p = 0 \}.
\]
The third analogue is a slight generalization of a construction given in \cite{HHexponent}.\footnote{Hochster and Huneke defined
the invariants $T_I$ and $T_I(x)$ for ideals $I \subseteq R$ and $x \in R$.  Their notions 
are in fact
a special case of what I define in this section.  In my notation, one would write 
$T_I^R$ and $T_I^R(x)$, considering $I$ to be a submodule of $R$.  However, all of the results
in \cite{HHexponent} for $T_I$ and $T_I(x)$ go through without change for
$T_N^M$ and $T_N^M(x)$ respectively.  It may be puzzling at first what to do with the hypothesis that 
$R$ is reduced.  However, in all cases one can eliminate it, sometimes with slight changes
in the statements.  For example, one must replace ``square locally stable test element'' everywhere
with ``$(q_0+1)$-power of a locally stable $q_0$-weak test element,'' and in Proposition 3.3(d)-(g), we must constrain $q$ to be $\geq q_0$.  I encourage the reader to read section 3 of \cite{HHexponent}
with these facts in mind.}  If $N$ is a submodule of $M$ and $z \in M$, then $\q \in \Spec R$ is
a \emph{stable prime associated to $N \subseteq M$ and $z$} if $z \notin (N_\q)^*_{M_\q}$ and for 
all primes $\p \subsetneq \q$, $z \in (N_\p)^*_{M_\p}$.  We write $\q \in T^M_N(z)$ in this case.
Set $T^M_N := \bigcup_{z \in M} T^M_N(z)$, the \emph{stable primes associated to $N \subseteq M$}.
Then our final tight closure theoretic notion of associated primes of $M$ is \[
T^M := T^M_0,
\]
the set of stable primes associated to $0 \subseteq M$.

How do we compare these three notions?  First, we have the following
\begin{prop}
If $R$ has a locally stable $q_0$-weak test element $c$, then $T^M \subseteq \left( \GAss{M} \right) \cap \Assph{M}.$
\end{prop}

\begin{proof}
Let $\q \in T^M$.  Then $\q \in T^M_0(z)$ for some $z \in M$.  In particular, $\q$
is minimal over the ideal $0^*_{F^{e'}(M)} : c z^{q'}$ for all $e' \gg 0$.  Hence,
$\q \in \Ass G^{e'}_R(M) \subseteq \GAss{M}$.

Now, suppose that $c z^q \in 0^*_{F^e(M_\q)}$ for all $q \gg 0$.  Then $c^{q_0 + 1} z^{q q_0} = 0$
for all $q \gg 0$, from which it follows that $z \in 0^*_{M_\q}$, contradicting the
definition of $T^M_0(z)$.  Hence, there are infinitely many powers $q$ of $p$ for which
$c z^q \notin 0^*_{F^e(M_\q)}$.  For all such $q$, then, $0^*_{F^e(M_\q)} : c z^q \subseteq \q R_\q$.
On the other hand, for any $\p \subsetneq \q$, $z \in 0^*_{M_\p}$, so 
$c z^q = 0$ in $M_\p = (M_\q)_{\p R_\q}$ for all $q \geq q_0$.  For these
choices of $q$, then, $\p R_\q \nsupseteq 0^*_{F^e(M_\q)} : c z^q$, so 
$\q R_\q$ is minimal over $0^*_{F^e(M_\q)} : c z^q$.  Therefore, $\q R_\q
\in \ffunc{Ass}_{R_\q} G^e_{R_\q}(M_\q)$, which proves that $\pdp{R_\q}
M_\q = 0$, i.e. $\q \in \Assph M$.
\end{proof}

For further comparisons, we examine the notion of phantom assassinator more 
closely:

\subsection*{Phantom assassinators}
Many properties of phantom assassinators were derived and used in \cite{AHH}, although
as mentioned before, it was assumed that the modules investigated had finite phantom
projective dimension.  It is natural to ask which of their results still hold in our
more general case involving modules of possibly infinite phantom projective dimension:

First of all, the
equivalence of (a): (that $\pdp{\m} M = 0$) and (c): (that $\m \in \cup_e \func{Ass} 
G^e(M)$) for a finitely generated module $M$ over a local ring $(R,\m)$ in Proposition 5.9
of that paper has already been shown here in Lemma~\ref{lem:phass},
assuming that $R$ has a weak test element.  Part (b) (that
$\ffunc{ppd}_R M = \func{mnht} \m$)
is of course false in the case of infinite phantom projective dimension, since the 
minheight of an ideal in a Noetherian ring is always finite,
bounded as it is by the height, hence by the size of a generating set of the ideal

Next we investigate \cite[Theorem 5.11]{AHH}.  The given proof of 
part (a) involves phantom resolutions, but in fact the first two statements of (a) can
be proved without them, as we shall show by using the following easy
\begin{lemma}\label{lem:Geincrease}\cite[Lemma 1.3]{nmepdepBC}
Let $(R,\m)$ be a Noetherian local ring of characteristic $p>0$ and $M$ a finitely
generated $R$-module.  Then for any $e \geq 0$, $\func{Ass} G^e(M) \subseteq
\func{Ass} G^{e+1}(M)$.
\end{lemma}

\begin{lemma}[Most of {\cite[Theorem 5.11(a)]{AHH}}, without assuming finite ppd or acceptable type]
Let $M$ be a finitely generated $R$-module, where $R$ is a 
Noetherian local ring of 
characteristic $p>0$ containing a locally stable weak test element.  Every minimal prime of $M$
is in $\Assph M$.  For all $e \in \N$, $\Assph F^e(M) = \Assph M$.
\end{lemma}
\begin{proof}
For the first statement, first assume that $\m$ is itself a minimal prime of $M$.
Then for some positive integer $t$, $\m^t M = 0$.  Take any $z \in M \setminus 0^*_M$.
Then $\m^t z = 0$, so $\m^t \subseteq 0^*_M : z \subseteq \m$.  The set in the middle
is the annihilator of $\bar{z}$ in $G^0(M)$, so it must be contained in some associated
prime of $G^0(M)$.  But the only prime ideal containing $\m^t$ is $\m$.  Hence
$\m \in \func{Ass} G^0(M)$.  By Lemma~\ref{lem:phass}, then, $\m$ consists of phantom zerodivisors
of $M$, so that $\pdp{\m} M = 0$, i.e. $\m \in \Assph M$.

Now take the general case where $\p$ is a minimal prime of $M$.  Then $\p R_\p$ is
a minimal prime of $M_\p$, and so $\pdp{\p R_\p} M_\p = 0$, which is the same as
saying that $\p \in \Assph M$.

For the second statement, note that for any $e$, Lemma~\ref{lem:Geincrease}
implies that \[
\bigcup_{e' \geq 0} \ffunc{Ass}_{R_\p} G^{e'}(M_\p) = \bigcup_{e' \geq e} 
\ffunc{Ass}_{R_\p} G^{e'}(M_\p).
\]
Moreover, $G^{e'}(F^e(M_\p)) = G^{e+e'}(M_\p)$ for any $e, e' \geq 0$, so by
Lemma~\ref{lem:phass}, the union of the left hand side of the equation above consists of 
the phantom zerodivisors of $M_\p$, whereas the union of the right hand side
consists of the phantom zerodivisors of $F^e(M_\p)$.  Thus, $M_\p$ has phantom
depth 0 if and only if $F^e(M_\p) = (F^e(M))_\p$ has phantom depth 0, which
shows that $\Assph M = \Assph F^e(M)$ for any $e \geq 0$.
\end{proof}

Part (b) of \cite[Theorem 5.11]{AHH} states
that if $M$ has finite phantom
projective dimension, then $\Assph M$ is a finite set.  In particular,
every element of $\Assph M$ is a minimal prime of some ideal of the form $I_i + \p_j$, where
$I_i = I_{r_i}(G.)$ for a given finite 
phantom resolution $G.$ of $M$ and $\p_j$ is a minimal prime of $R$.  In general,
there will at least be non-local cases where $\Assph M$ is an 
infinite set (see below), so I doubt that \cite[Theorem 5.11(b)]{AHH} has
any analogue at all in the case of infinite ppd.

Now consider \cite[Theorems 5.13-5.14]{AHH}.  The only results used in the proofs are
ones which go through without assuming $M$ to have finite ppd, as the reader may check.  So we state the
results here (with the second one in slightly different form, for the sake of
simplicity):

\begin{prop}[{\cite[Theorem 5.13]{AHH}}, without assuming finite ppd or acceptable type]
\label{prop:Gassphass}
Let $R$ be a Noetherian local ring of characteristic $p$ possessing a locally stable weak
test element, and let $M$ be a finitely generated $R$-module.  Then every 
maximal element of $\bigcup_{e \in \N} \func{Ass} G^e(M)$ is in the phantom
assassinator of $M$.
\end{prop}

By the Singh-Swanson example mentioned in note~\ref{av}, $\bigcup_{e \in \N} \func{Ass}
G^e(M)$ may have infinitely many maximal elements, even for very nice rings $R$.  Singh and
Swanson's ring is not
local, but if we could make it local, we would have an example of a finitely generated
$M$ over an otherwise well-behaved local ring where $\Assph M$ was an infinite set.

\begin{definition}[See {\cite[5.8]{AHH}}]
We say that $M$ is \emph{phantom unmixed} if every prime in $\Assph M$ is
a minimal prime of $M$.
\end{definition}

\begin{prop}[{\cite[Theorem 5.14]{AHH}}, without assuming finite ppd or acceptable type]\label{prop:phunmix}
Let $R$ be a Noetherian local ring of characteristic $p>0$ possessing
a locally stable weak test element, and let $M$ be a finitely generated $R$-module.
Then $M$ is phantom unmixed if and only if $G^e(M)$ is unmixed
for all $e \geq 0$.

Hence (because of \cite[Theorem 3.7]{AHH}), if $M$ is phantom unmixed
then tight closure commutes with localization for the pair $0
\subseteq M$.
\end{prop}

This last result is not new.  As Craig Huneke pointed out to me, it is a consequence of
the following theorem:

\begin{thm}(\cite[Theorem 3.5]{HHexponent}, in terms of modules)
Let $R$ be a Noetherian local ring, let $c$ be the $(q_0 + 1)$'th power of a
locally stable $q_0$-weak test element, and let $M$ be a finitely generated $R$-module.
The following are equivalent: \begin{enumerate}
\item Tight closure commutes with localization for the pair $0 \subseteq M$.
\item \begin{enumerate}
 \item $T^M$ is finite and 
 \item For every $z \in M$ and $\q \in T^M(z)$, there exists an integer $N$, possibly 
 depending on $\q$, such that for all $q \gg 0$,
   $\q^{N q} \subseteq \left(0^*_{F^e(M)} : c z^q \right)_\q.$
\end{enumerate}
\end{enumerate}
\end{thm}

In the case of a phantom unmixed module $M$, the set $\Assph M$ is finite, so the smaller
set $T^M$ is certainly finite.  Also, since any $\q \in \Assph M$  is a minimal prime of $M$,
the same is true for any $\q \in T^M$,
so for any such $\q$ there is some positive integer $L$ with $\q^L M_\q = 0$.  Let $N = L + \mu(\q) - 1$.  Then
$\q^{N q} \subseteq (\q^{[q]})^{N - \mu(\q) + 1} = (\q^L)^{[q]}$, so that 
$\q^{N q} F^e(M_\q) \subseteq (\q^L M_\q)^{[q]}_{M_\q} = 0$ for any power $q$ of $p$, which
certainly proves (b).  Then by the theorem, tight closure commutes with localization for
$0 \subseteq M$.

\section{Diagram-chasing with exponents}
\label{sec:pseudolesAPX}
Given a short exact sequence of complexes of modules, classical homological algebra shows us that
there is an associated long exact sequence of homology modules, from which we can derive useful
criteria for vanishing and other properties of the homology modules in question.  On the other hand,
if we only have a short stably phantom exact sequence of complexes of finitely-generated $R$-modules,
I do not know how to obtain a useful long sequence of modules.  However, in some important ways,
we may act as if we had such a long sequence, as I make precise in the main result of this
section, Proposition~\ref{prop:pseudolesAPX}.  Along the way we will define ``phantom connecting 
homomorphisms,'' and I hope that these along with Proposition~\ref{prop:pseudolesAPX} will be 
useful additions to the ever-expanding toolbox of tight closure theory.

Throughout this section, let $R$ be a Noetherian ring of prime characteristic $p>0$ which contains a
$q_0$-weak test element $c$, and let \[
0 \rightarrow L. \stackrel{\alpha.}{\rightarrow} M. \stackrel{\beta.}{\rightarrow} N. \rightarrow 0
\]
be a sequence of complexes of finitely generated $R$-modules, which is a complex
in each degree, and such that for
each $i \in \Z$ and every positive integer $e \geq 0$, we have: \begin{itemize}
\item $\beta_i$ is surjective,
\item $\ker F^e(\beta_i) \subseteq \left( \im F^e(\alpha_i) \right)^*_{F^e(M_i)}$, and
\item $\ker F^e(\alpha_i) \subseteq 0^*_{F^e(L_i)}$.
\end{itemize}
In other words, we have a stably phantom exact complex in each degree.  We call this a short
stably phantom exact sequence of complexes.

The point is that although there seems to be no long exact sequence which arises
from this situation, we can get the same kind of benefits that we would get out of a
long exact sequence, ``up to tight closure.''

To simplify notation, we will abuse it by assuming the
reader can keep track of the
Frobenius powers on the maps, in the same manner as in Section~\ref{sec:pseudolesABBR}.

\subsection{Phantom connecting homomorphisms}
Fix an integer $i$ and $e \geq 0$, and let $[z] \in H_i(F^e(N.))$.  That is,
$d_i^N(z) = 0.$  Since $\beta_i$ is surjective, there is some $y \in F^e(M_i)$
with $\beta_i(y) = z$.  Now, \[
\beta_{i-1}(d_i^M(y)) = d_i^N(\beta_i(y)) = d_i^N(z) = 0,
\]
so $d_i^M(y) \in \ker \beta_{i-1} \subseteq (\im \alpha_{i-1})^*_{F^e(M_{i-1})}$.
Therefore, for any $e' \geq e_0$, we have $d_i^M(c y^{q'}) = c \cdot d_i^M(y)^{q'} 
\in \im \alpha_{i-1}$.  Fixing $e'$, this means that there is some $x \in F^{e+e'}
(L_{i-1})$ such that $\alpha_{i-1}(x) = d_i^M(c y^{q'})$.  When $x$ arises from
$[z]$ in this way, we write: \[
x \in \delta_i^{(q')}([z]).
\]
That is, \[
\delta_i^{(q')}([z]) = \{ x \in F^{e+e'}(L_{i-1}) \mid \exists y \in F^e(M_i):
\beta_i(y) = z \text{ and } d_i^M(c y^{q'}) = \alpha_{i-1}(x) \}.
\]

At this point it would be tempting to say that our ``connecting homomorphism''
sends, for each $q'$, $[z]$ to an $[x]$ which arises in the above manner.  However,
it is not even clear that $d_{i-1}^L(x) = 0$, so we might not even get into the 
$(i-1)$'th homology module of $F^{e+e'}(L.)$.  Even if we do, it is not clear
that the homology class is unique.  However, both of these things are true ``up
to tight closure''.  More precisely, 
\begin{prop}
Let $z \in H_i(F^e(N.))$, and let $e', e'' \geq e_0$.  Then for
any $x \in \delta_i^{(q')}([z])$, $d_{i-1}^L(c x^{q''}) = 0$.  Moreover, the homology class
$[c x^{q''}] \in H_{i-1}(F^{e+e'+e''}(L.))$ is independent of the choice of $x$.
\end{prop}

\begin{proof}
We have \[
\alpha_{i-2}(d_{i-1}^L(x)) = d_{i-1}^M(\alpha_{i-1}(x)) = 
d_{i-1}^M(d_i^M(c y^{q'})) = 0,
\]
and thus, $d_{i-1}^L(x) \in \ker \alpha_{i-2} \subseteq 0^*_{F^{e+e'}(L_{i-2})}$.
Hence, for any fixed $e'' \geq e_0$, $d_{i-1}^L(c x^{q''}) = c \cdot d_{i-1}^L(x)^{q''}
= 0$.  As for uniqueness
up to tight closure, let $y_1 \in F^e(M_i)$ and $x_1 \in F^{e+e'}(L_{i-1})$
such that $\beta_i(y_1) = z$ and $\alpha_{i-1}(x_1) = d_i^M(c y_1^{q'})$.  Then
$\beta_i(y - y_1) = 0$, so $y-y_1 \in (\im \alpha_i)^*_{F^{e+e'}(M_i)}$, whence 
$c y^{q'} - c y_1^{q'} 
\in \im \alpha_i$.  Say $c y^{q'} - c y_1^{q'} = \alpha_i(w)$.  Then \[
\alpha_{i-1}(x - x_1) = d_i^M(c y^{q'} - c y_1^{q'}) = d_i^M(\alpha_i(w)) = 
\alpha_{i-1}(d_i^L(w)).
\]
Thus, $x - x_1 - d_i^L(w) \in \ker \alpha_{i-1} \subseteq 0^*_{F^{e+e'}(L_{i-1})}$, 
so that $c x^{q''} - c x_1^{q''} = d_i^L (c w^{q''})$, and hence 
$\left[c x^{q''}\right] = \left[c x_1^{q''}\right]$.
\end{proof}

With this existence and uniqueness in mind, for any $e', e'' \geq e_0$ we set \[
\delta_i^{(q',q'')}([z]) = \left[c x^{q''}\right],
\]
where $x$ is an arbitrarily chosen element of $\delta_i^{(q')}([z])$.  This defines
an ${}^{e'+e''}R$-linear homomorphism of homology modules, and we call it a \emph{phantom
connecting homomorphism}.

The following lemma is useful in calculations:

\begin{lemma}\label{lem:deltaphzero}
Let $e \geq 0$, $e' \geq e_0$, and $[z] \in H_i(F^e(N.))$, and suppose that for
all $e'' \gg 0$, $\delta_i^{(q',q'')}([z])$ is a phantom element of 
$H_{i-1}(F^{e+e'+e''}(L.))$.  Then for all $e'' \geq e_0$, $\delta_i^{(q',q'')}([z]) 
= 0$.
\end{lemma}
\begin{proof}
Let $x \in \delta_i^{(q')}([z])$.  For all large $q''$, 
$\delta_i^{(q',q'')}([z]) = [c x^{q''}]$
is phantom.  That is, $c x^{q''} \in (\im d_i^L)^*_{F^{e+e'+e''}(L_{i-1})}$.
In particular, then, $c^{q_0+1} x^{q''q_0} \in \im d_i^L$.  It follows that 
$x \in (\im d_i^L)^*_{F^{e+e'}(L_{i-1})}$.
Therefore, for all $e'' \geq e_0$, $c x^{q''} \in \im d_i^L$,
so that $\delta_i^{(q',q'')}([z]) = [c x^{q''}] = 0.$
\end{proof}

\subsection{Diagram-chasing with exponents}

We start with a lemma characterizing elements that get sent to phantom elements:

\begin{lemma}[What gets sent to phantom elements]\label{lem:pseudoles1}
In the setup of this section, we have the following for every $i \in \Z$:
\begin{enumerate}[(a)]
\item\label{it:Lexact} Let $[x] \in H_i(F^e(L.))$.  Then
$H_i(F^e(\alpha.))([x])$ is a phantom element of $H_i(F^e(M.))$
if and only if for every $e' \geq e_0$, there exists $[z] \in H_{i+1}(F^{e+e'}(N.))$ such that $c_1 x^{q' q_0} \in \delta_{i+1}^{(q_0)}([z]).$
\item\label{it:Mexact} Let $[y] \in H_i(F^e(M.))$.  Then $H_i(F^e(\beta.))([y])$ is a phantom
element of $H_i(F^e(N.))$ if and only if for all $e' \geq 3 e_0$, 
$[c_2 y^{q'}] \in \im H_i(F^{e+e'}(\alpha.))$.
\item\label{it:Nexact} Let $[z] \in H_i(F^e(N.))$, and fix $e' \geq e_0$.  
Then $\delta_i^{(q',q'')}([z]) =0$ for all $e'' \geq e_0$ if and only if 
$[c^{q''+1} z^{q'q''}] \in \im H_i(F^{e+e'+e''}(\beta.))$
for all $e'' \geq e_0$.
\end{enumerate}
\end{lemma}

\noindent \emph{Remark:} This Lemma is an analogue
of the fact that in a three-term exact
sequence, the image of the first map equals the kernel of the second map.  

\noindent \emph{Remark:} In practice, $c$ may often be taken to be a test element and not merely a weak test element.  In such cases, in the above lemma and the following Proposition one may replace every occurrence of $e_0$ with $0$ and every occurrence of $c_n$ with $c^{n+1}$.  For example, if $c$ is a test element, part (b) says that $H_i(F^e(\beta.))([y])$ is phantom if and only if for all $e' \geq 0$, $c^3 y^{q'} \in \im H_i(F^{e+e'}(\alpha.))$.

\begin{proof}[Proof (\ref{it:Lexact}).]
First suppose that $[\alpha_i(x)]$ is a phantom element of $H_i(F^e(M.))$.  This means that
for all $e' \geq e_0$, there exists $y \in F^{e+e'}(M_{i+1})$ with $\alpha_i(c x^{q'}) =
d_{i+1}^M(y)$.  Also, \[
d_{i+1}^N(\beta_{i+1}(y)) = \beta_i(d_{i+1}^M(y)) = \beta_i(\alpha_i(c x^{q'})) = 0,
\]
so $[\beta_{i+1}(y)] \in H_{i+1}(F^{e+e'}(N.))$.  It follows 
that $c (c x^{q'})^{q_0} = c_1 x^{q'q_0} \in \delta_{i+1}^{(q_0)}([\beta_{i+1}(y)])$.

Conversely suppose that for any $e' \geq e_0$, there exists $[z] \in 
H_{i+1}(F^{e+e'}(N.))$ such that $c_1 x^{q'q_0} \in \delta_{i+1}^{(q_0)}([z]).$
By definition, then, for each such $e'$ there exists some $y \in F^{e+e'}(M_{i+1})$ such that $\beta_{i+1}(y) = z$
and \[
c_1 \alpha_i(x)^{q' q_0} = \alpha_i(c_1 x^{q' q_0}) = d_{i+1}^M(c y^{q_0}) \in (\im d_{i+1}^M)^{[q' q_0]}_{F^e(M)}
\]  Since this holds for any $e' \geq e_0$, 
it follows that $\alpha_i(x) \in (\im d_{i+1}^M)^*_{F^e(M_i)}$.
That is, $[\alpha_i(x)]$ is a phantom element of $H_i(F^e(M.))$.
\end{proof}
\begin{proof}[Proof (\ref{it:Mexact}).]
First suppose that $H_i(F^e(\beta.))([y])$ is a phantom element of $H_i(F^e(N.))$.  Then
for any $e'' \geq e_0$, since $\beta_{i+1}$ is surjective there is some 
$w \in F^{e+e''}(M_{i+1})$ with $\beta_i(c y^{q''})
= d_{i+1}^N(\beta_{i+1}(w))$.  Then we have: \[
\beta_i(c y^{q''}) = d_{i+1}^N(\beta_{i+1}(w)) = \beta_i(d_{i+1}^M(w)),
\]
so $c y^{q''} - d_{i+1}^M(w) \in \ker \beta_i \subseteq (\im \alpha_i)^*_{F^{e+e''}(M_i)}$.
Hence, there is some $x \in F^{e+e''+e_0}(L_i)$
with $\alpha_i(x) = c^{q_0+1} y^{q''q_0} - d_{i+1}^M(c w^{q_0})$.  We have
\begin{align*}
\alpha_{i-1}(d_i^L(x)) &= d_i^M(\alpha_i(x)) = d_i^M(c^{q_0+1} y^{q''q_0} - d_{i+1}^M
(c w^{q_0})) \\
&= c^{q_0+1} (d_i^M(y))^{q''q_0} - d_i^M(d_{i+1}^M(c w^{q_0})) = 0.
\end{align*}
Thus, $d_i^L(x) \in \ker \alpha_{i-1} \subseteq 0^*_{F^{e+e''+e_0}(L_{i-1})}$.  It follows
that $d_i^L(c x^{q_0}) = 0$.  Let $e' = e'' + 2e_0$.  Then $[c x^{q_0}] \in 
H_i(F^{e + e'}(L.))$, and \[
H_i(F^{e+e'}(\alpha.))([c x^{q_0}]) = [c_2 y^{q''q_0^2} - 
d_{i+1}^M(c^{q_0 + 1}w^{q_0^2})] = [c_2 y^{q'}],
\]
which proves the forward direction of part~(\ref{it:Mexact}).

Conversely, suppose that for all $e' \geq 3 e_0$ there exists 
$[x] \in H_i(F^{e+e'}(L.))$
with $H_i(F^{e+e'}(\alpha.))([x]) = [c_2 y^{q'}].$  Then for each such $q'$,
there exists $w \in F^{e+e'}(M_{i+1})$ with \[
c_2 y^{q'} = \alpha_i(x) + d_{i+1}^M(w).
\]
Then applying $\beta_i$ to both sides, we have \[
c_2 \beta_i(y)^{q'} = d_{i+1}^N(\beta_{i+1}(w)) \in \im F^{e+e'}(d_{i+1}^N)
= (\im F^e(d_{i+1}^N))^{[q']}_{F^e(N_i)},
\]
and since this holds for all sufficiently large powers $q'$ of $p$, it follows
that $[\beta_i(y)]$ is a phantom element of
$H_i(F^e(N.))$.
\end{proof}

\begin{proof}[Proof (\ref{it:Nexact})]
First, fixing any $e' \geq e_0$ and $e'' \geq e_0$, suppose that 
$\delta_i^{(q',q'')}([z]) = 0$.  Then there are $y$, $x$,
and $v$ such that $\beta_i(y) = z$, $\alpha_{i-1}(x) = d_i^M(c y^{q'})$,
and $c x^{q''} = d_i^L(v)$.  Then \begin{align*}
d_i^M(\alpha_i(v)) &= \alpha_{i-1}(d_i^L(v)) = \alpha_{i-1}(c x^{q''}) 
= c \alpha_{i-1}(x)^{q''}\\
 &= c d_i^M(c y^{q'})^{q''} = d_i^M(c^{q''+1} y^{q'q''}).
\end{align*}
Let $y'' = c^{q''+1} y^{q'q''} - \alpha_i(v)$.  Then $d_i^M(y'') = 0$
and $\beta_i(y'') = c^{q''+1} z^{q'q''}$.  That is, \[
[c^{q''+1} z^{q'q''}] = H_i(F^{e+e'+e''}(\beta.))([y'']).
\]

For the converse, suppose that for all $e'' \geq e_0$, $[c^{q''+1} z^{q'q''}]
\in \im H_i(F^{e+e'+e''}(\beta.))$.  Then fixing $e''$, there exist $a \in 
F^{e+e'+e''}(M_i)$ and $b \in F^{e+e'+e''}(N_{i+1})$ such that $d_i^M(a) = 0$ and
$c^{q''+1} z^{q'q''} = \beta_i(a) + d_{i+1}^N(b)$.  However, since $\beta_i$ 
and $\beta_{i+1}$ are surjective, there exist $g$ and $y$ such that 
$z = \beta_i(y)$ and $b = \beta_{i+1}(g)$.  Then if we set $v = a + d_{i+1}^M(g)$,
it follows that $d_i^M(v) = 0$ and \[
\beta_i(c^{q''+1} y^{q'q''} - v) = c^{q''+1} z^{q'q''} - \beta_i(a) - 
d_{i+1}^N(\beta_{i+1}(g)) = 0.
\]
Hence, $c^{q''+1} y^{q'q''} - v \in \ker \beta_i \subseteq (\im 
\alpha_i)^*_{F^{e+e'+e''}(M_i)}$.  In particular, there is some 
$t \in F^{e+e'+e''+e_0}(L_i)$ with \[
\alpha_i(t) = c^{q''q_0 + q_0 + 1} y^{q'q''q_0} - c v^{q_0}.
\]
Now, by definition, $\delta_i^{(q',q'')}([z]) = [c x^{q''}]$, where
$\alpha_{i-1}(x) = d_i^M(c y^{q'})$.  Then \begin{align*}
\alpha_{i-1}(c^{q_0+1} x^{q''q_0}) &= c^{q_0+1} d_i^M(c y^{q'})^{q''q_0}
= d_i^M(c^{q''q_0 + q_0 + 1} y^{q'q''q_0} - c v^{q_0}) \\
&= d_i^M(\alpha_i(t)) = \alpha_{i-1}(d_i^L(t)).
\end{align*}
Thus, $c^{q_0+1} x^{q''q_0} - d_i^L(t) \in \ker \alpha_{i-1} \subseteq 0^*$,
whence $c_2 x^{q''q_0^2} = d_i^L (c t^{q_0})$.  Since this
holds for all $q'' \geq q_0$, $x \in (\im d_i^L)^*_{F^e(L_{i-1})}$.  Thus, 
for all $q'' \geq q_0$, $c x^{q''} \in \im d_i^L$.  That is, \[
\delta_i^{(q',q'')}([z]) = [c x^{q''}] = 0.
\]
\end{proof}

Now we come to the main result of this Section:

\begin{prop}[The spots where stably phantom homology occurs]\label{prop:pseudolesAPX}
In the setup of this section, we have the following for each $i \in \Z$: 
\begin{enumerate}[(a)]
\item\label{it:Lph} $L.$ has stably phantom homology at $i$ if and only if the
following two conditions hold for all $e \geq 0$:
\begin{enumerate}[i.]
  \item Any element of $H_i(F^e(M.))$ which $H_i(F^e(\beta.))$ maps to a 
  phantom element of $H_i(F^e(N.))$ is itself a phantom element of $H_i(F^e(M.))$,
  and
  \item For any $[z] \in H_{i+1}(F^e(N.))$ and any integers $e' \geq e_0$ and
  $e'' \geq e_0$,
  $[c^{q'' + 1} z^{q'q''}] \in \im H_{i+1}(F^{e+e'+e''}(\beta.))$.
\end{enumerate}
\item\label{it:Mph} $M.$ has stably phantom homology at $i$ if and only if the
following two conditions hold for all $e \geq 0$:
\begin{enumerate}[i.]
  \item Any $[z] \in H_i(F^e(N.))$ for which $\delta_i^{(q',q'')}([z]) = 0$ for 
  all $e' \geq e_0$ and $e'' \geq e_0$ is a phantom element of $H_i(F^e(N.))$, and
  \item For any $[x] \in H_i(F^e(L.))$ and any $e' \geq e_0$ 
  there exists $[z] \in H_{i+1}(F^{e+e'}(N.))$ such that $c_1 x^{q'q_0} \in
  \delta_{i+1}^{(q_0)}([z])$.
\end{enumerate}
\item\label{it:NphAPX} $N.$ has stably phantom homology at $i$ if and only if the
following two conditions hold for all $e \geq 0$:
\begin{enumerate}[i.]
  \item Any element of $H_{i-1}(F^e(L.))$ which $H_{i-1}(F^e(\alpha.))$ maps
  to a phantom element of $H_{i-1}(F^e(M.))$ is itself a phantom element of 
  $H_{i-1}(F^e(L.))$, and
  \item For any $[y] \in H_i(F^e(M.))$ and any integer $e' \geq 3 e_0$,
  $[c_2 y^{q'}] \in \im H_i(F^{e+e'})(\alpha.)$.
\end{enumerate}
\item\label{it:MNph} If $N.$ has stably phantom homology at $i+1$ and $M.$ has stably 
phantom homology at $i$, then 
$L.$ has stably phantom homology at $i$.
\item\label{it:LNph} If both $L.$ and $N.$ have stably phantom homology at $i$, so does $M.$.
\item\label{it:LMphAPX} If $M.$ has stably phantom homology at $i$ and $L.$ has stably phantom 
homology at $i-1$, then $N.$ has stably phantom homology at $i$.
\end{enumerate}
\end{prop}

Parts $a$, $b$, and $c$ 
are analogues of the fact that in a five-term exact sequence, the middle term is zero if and 
only if the first map is surjective and the last map is injective.  Note that part $c$ is 
essentially a generalization of Proposition~\ref{prop:pseudolesABBR}(1).

Parts $d$, $e$, and $f$ are analogues
of the fact that in a three-term exact sequence, if the outer two terms vanish, so does the
middle term.  Note that part $f$ is a generalization of Proposition~\ref{prop:pseudolesABBR}(2).

\begin{proof}[Proof (\ref{it:Lph})]
First suppose that $L.$ has stably phantom homology at $i$.

To prove 
condition (i), let $[y] \in H_i(F^e(M.))$ such that $[\beta_i(y)]$ 
is phantom.  Then by Lemma~\ref{lem:pseudoles1}(\ref{it:Mexact}), for every $e' \geq 3e_0$ there exists $[x] \in 
H_i(F^{e+e'}(L.))$ such that $[c_2 y^{q'}] = [\alpha_i(x)]$.
But $[x]$ is phantom, so $[c x^{q_0}] = 0$, which implies that
$[c_3 y^{q'q_0}] = [\alpha_i(c x^{q_0})] = 0.$  Since this
holds for all $e' \geq 3 e_0$, $[y]$ is phantom.

To prove condition (ii), let $[z] \in H_{i+1}(F^e(N.))$.  
Then for all $e'' \geq e_0$, 
$\delta_i^{(q',q'')}([z])$ is phantom, hence by Lemma~\ref{lem:deltaphzero},
$\delta_i^{(q',q'')}([z]) = 0$ for all $e'' \geq e_0$.  Then by 
Lemma~\ref{lem:pseudoles1}(\ref{it:Nexact}), $[c^{q''+1} z^{q' q''}] \in 
\im H_{i+1}(F^{e+e'+e''}(\beta.))$ for all $e'' \geq e_0$.

Conversely, suppose that conditions (i) and (ii) hold for all $e \geq 0$.  Let
$[x] \in H_i(F^e(L.))$.  Then $H_i(F^e(\beta.))([\alpha_i(x)]) =
[\beta_i(\alpha_i(x))] = 0$, so it is certainly a phantom element of
$H_i(F^e(N.))$.  Thus, by condition (i), $[\alpha_i(x)]$ is a phantom
element of $H_i(F^e(M.))$.  Now by Lemma~\ref{lem:pseudoles1}(\ref{it:Lexact}), for any $e' \geq e_0$
there is some $[z] \in H_{i+1}(F^{e+e'}(N.))$ with 
$c^{q_0+1} x^{q'q_0} \in \delta_{i+1}^{(q_0)}([z])$.  In particular, for some
$y \in F^{e+e'}(M_{i+1})$, $\beta_{i+1}(y) = z$ and $d_{i+1}^M(c y^{q_0})
= \alpha_i(c^{q_0+1}x^{q'q_0})$.  Condition (ii) guarantees that there is some $[v] \in 
H_{i+1}(F^{e+e'+2e_0}(M.))$ with $[c^{q_0 + 1} z^{q_0^2}] = [\beta_{i+1}(v)]$.
Since $\beta_{i+2}$ is surjective, it follows without loss of generality
that $\beta_{i+1}(c^{q_0+1} y^{q_0^2} - v) = 0$.  That is, \[
c^{q_0+1} y^{q_0^2} - v \in \ker \beta_{i+1} \subseteq (\im 
\alpha_{i+1})^*_{F^{e+e'+2e_0}(M_{i+1})}.
\]
In particular, there is some $u$ with $\alpha_{i+1}(u) = c_2
y^{q_0^3} - c v^{q_0}$.  Then \begin{align*}
\alpha_i(d_{i+1}^L(u)) &= d_{i+1}^M(\alpha_{i+1}(u)) = d_{i+1}^M(c_2 y^{q_0^3} 
- c v^{q_0}) \\
&= c_1 d_{i+1}^M (c y^{q_0})^{q_0^2} = c_1 \alpha_i(
c_1 x^{q'q_0})^{q_0^2} = \alpha_i ( c_3 x^{q' q_0^3}).
\end{align*}
That is, $c_3 x^{q' q_0^3} - d_{i+1}^L(u) \in 
\ker \alpha_i \subseteq 0^*_{F^{e+e'+3e_0}(L_i)}$.  In particular,
$c_4 x^{q' q_0^4} = d_{i+1}^L(c u^{q_0}) \in \im d_{i+1}^L.$
Since this holds for all large $q'$, it follows that $x \in (\im 
d_{i+1}^L)^*_{F^e(L_i)}$, so that $[x]$ is a phantom element of
$H_i(F^e(L.))$.
\end{proof}

\begin{proof}[Proof (\ref{it:Mph})]
First suppose that $M.$ has stably phantom homology at $i$.

To prove condition (i), let $[z] \in H_i(F^e(N.))$ such that
for all $e', e'' \geq e_0$, $\delta_i^{(q',q'')}([z]) = 0$.  By 
Lemma~\ref{lem:pseudoles1}(\ref{it:Mexact}), $[c^{q''+1} z^{q'q''}] \in 
\im H_i(F^{e+e'+e''}(\beta.))$ for all $e',e'' \geq e_0$.  In particular,
there is some $[y] \in H_i(F^{e + e'+e_0}(M.))$ such that 
$[c^{q_0+1} z^{q'q_0}] = [\beta_i(y)]$.  Since $[y]$ is
phantom, $[c y^{q_0}] = 0$, which implies that 
$\left[c_2 z^{q'q_0^2}\right] = [\beta_i(c y^{q_0})] = 0$.  Since this holds for all
$q' \geq q_0$, it follows that $[z]$ is a phantom element of $H_i(F^e(N.))$.

To prove condition (ii), note that for any $[x] \in H_{i-1}(F^e(L.))$,
$[\alpha_i(x)]$ is phantom.  Then the conclusion follows directly from
Lemma~\ref{lem:pseudoles1}(\ref{it:Lexact}).

Conversely, suppose that conditions (i) and (ii) hold for all $e \geq 0$.  Let
$[y] \in H_i(F^e(M.))$.  For all $e',e'' \geq e_0$, we have \[
[c^{q''+1} \beta_i(y)^{q'q''}] = H_i(F^{e+e'+e''}(\beta.))([c^{q''+1} y^{q'q''}]),
\]
so that by Lemma~\ref{lem:pseudoles1}(\ref{it:Nexact}), $\delta_i^{(q',q'')}([\beta_i(y)])=0$
for all $e',e'' \geq e_0$.  Condition (i) then shows us that $[\beta_i(y)]$ is
a phantom element of $H_i(F^e(N.))$.  Therefore by Lemma~\ref{lem:pseudoles1}(\ref{it:Mexact}), 
$[c_2 y^{q'}] \in \im H_i(F^{e+e'}(\alpha.))$ for all $e' \geq 3 e_0$.
Temporarily fixing $e'$, say $[c_2 y^{q'}] = [\alpha_i(x)]$, where $d_i^L(x) = 0$.  Then by condition (ii),
for all $e'' \geq e_0$, there exists $[w] \in H_{i+1}(F^{e+e'+e''}(N.))$
with $c_1 x^{q'' q_0} \in \delta_{i+1}^{(q_0)}([w])$.  Now by
Lemma~\ref{lem:pseudoles1}(\ref{it:Lexact}), $[\alpha_i(x)] = [c_2 y^{q'}]$ is a phantom element of
$H_i(F^{e+e'}(M.))$  In particular, $c_3 y^{q' q_0} \in \im d_{i+1}^M$, and
since this holds for all large $q'$, it follows that $[y]$ is a phantom element
of $H_i(F^e(M.))$.
\end{proof}

\begin{proof}[Proof (\ref{it:NphAPX})]
First suppose that $N.$ has stably phantom homology at $i$.

To prove condition (i), let $[x] \in H_{i-1}(F^e(L.))$ be such that 
$[\alpha_{i-1}(x)]$ is a phantom element of $H_{i-1}(F^e(M.))$.  Then by
Lemma~\ref{lem:pseudoles1}(\ref{it:Lexact}), for every $e'\geq e_0$ there exists $[z] \in 
H_i(F^{e+e'}(N.))$ with $c_1 x^{q'q_0} \in \delta_i^{(q_0)}([z])$.
That is, there exists $y$ such that $\beta_i(y) = z$ and \[
d_i(c y^{q'}) = \alpha_{i-1}(c^{q'+1} x^{q'q''}).
\]
Since $[z]$ is phantom and $\beta.$ is surjective, there
exists $v$ such that $\beta_i(c y^{q_0}) = cz^{q_0} = d_i^N(\beta_{i+1}(v))
= \beta_{i-1}(d_i^M(v)).$  Thus, $\beta_i(c y^{q_0} - d_i^M(v)) = 0$, so that
there is some $t$ such that $c_1 y^{q_0^2} - d_i^M(c v^{q_0}) = \alpha_i(t)$.

Then \[
\alpha_{i-1}(d_i^L(t)) = c_1 d_i^M(y)^{q_0^2} = c d_i^M (c y^{q_0})^{q_0} 
= c \alpha_{i-1}(c_1 x^{q' q_0})^{q_0} = \alpha_{i-1}
(c_2 x^{q' q_0^2}).
\]
Hence, $c_2 x^{q' q_0^2} - d t \in \ker \alpha_i \subseteq 0^*$,
so $c_3 x^{q' q_0^3} \in \im d_i^L$.  Since $q'$ may be 
arbitrarily large, it follows that $[x]$ is phantom.

As for condition (ii), note that for any $[y] \in H_i(F^e(M.))$, $[\beta_i(y)]$ is 
phantom. Then the conclusion follows directly from Lemma~\ref{lem:pseudoles1}(\ref{it:Mexact}).

Conversely, suppose that conditions (i) and (ii) hold for all $e \geq 0$, and
let $[z] \in H_i(F^e(N.))$.  Then for some $y \in F^e(M_i)$, $\beta_i(y) = z$ and $\beta_{i-1}(d_i(y)) = d_i(z) = 0$, so $d_i(y) \in \ker \beta_{i-1} \subseteq (\im \alpha_{i-1})^*_{F^e(M_{i-1})}$.  Thus, for any fixed $e'\geq e_0$, there is some $x \in F^e(L_{i-1})$
with $d_i(c y^{q'}) = \alpha_{i-1}(x)$.

We have $\alpha_{i-2}(d_{i-1}^L(x)) = d_{i-1}^M(\alpha_{i-1}(x)) = 
d_{i-1}(d_i(c y^{q'})) = 0$, so $dx \in \ker \alpha_{i-2} \subseteq 0^*_{F^{e+e'}(L_{i-2})}$.  This means that for any $e'' \geq e_0$, $d_{i-1}^L(c x^{q''}) = 0$, so that $[c x^{q''}] \in H_{i-1}(F^{e+e'+e''}(L.))$.  On the other hand, since $\alpha_{i-1}(x) \in \im d_i^M$, it follows that $\alpha_{i-1}(c x^{q''}) \in \im d_i^M$, so that $[\alpha_{i-1}(c x^{q''})] = 0$.  In particular, \linebreak $H_{i-1}(F^{e+e'+e''}(\alpha.))([c x^{q''}]) = 0$, hence is phantom, so by condition (i), $[c x^{q''}]$ is a phantom element of $H_{i-1}(F^{e+e'+e''}(L.))$  Thus, $c_1 x^{q'' q_0} = c (c x^{q''})^{q_0} \in \im d_i^L$.  Since this holds for all large $q''$, it follows that $x \in (\im d_i^L)^*_{F^{e+e'}(L_i)}$.

In particular, there is some $w \in F^{e+e'}(L_{i+1})$ such that $c x^{q_0} = d_i^L(w)$.  We have \[
d \alpha(w) = \alpha d (w) = \alpha(c x^{q_0}) = d(c_1 y^{q' q_0}),
\]
so that $v := c_1 y^{q' q_0} - \alpha_i(w) \in \ker d_i^M$, whence
$[v] \in H_i(F^{e+e'+e_0}(M.))$.  By condition (ii), we get that $[c_2 v^{q_0^3}] \in \im H_i(F^{e + e' + 4 e_0}(\alpha.))$.  Say $c_2 v^{q_0^3} = \alpha_i(u) + d_{i+1}^M(t)$, where $d_i^M(u) = 0$.  Since $\beta_i(v) = c_1 z^{q' q_0}$, we have $\beta_i(c_2 v^{q_0^3}) = c_4 z^{q' q_0^4}$, whence we have \[
c_4 z^{q'q_0^4} = \beta_i(c_2 v^{q_0^3}) = \beta_i(\alpha_i(u)) + \beta_i(d_{i+1}^M(t)) \in \im d_{i+1}^N.
\]
Since this holds for all $e' \geq 0$, it follows that $z \in (\im d_{i+1}^N)^*_{F^e(N_i)}$, which is to say that $[z]$ is a phantom element of $H_i(F^e(N.))$.
\end{proof}

\begin{proof}[Proof (\ref{it:MNph}).]
Let $[x] \in H_i(F^e(L.))$.  Since every element of $H_i(F^e(M.))$ is phantom,
$H_i(F^e(\alpha.))([x])$ is phantom, so by part~(\ref{it:NphAPX})i, $[x]$ must be
a phantom element of $H_i(F^e(L.))$.
\end{proof}

\begin{proof}[Proof (\ref{it:LNph}).]
Let $[y] \in H_i(F^e(M.))$.  Then $H_i(F^e(\beta.))$ maps $[y]$ to a phantom element
of $H_i(F^e(N.))$, since all elements of the latter homology module are phantom, and 
thus by part~(\ref{it:Lph})i, $[y]$ is a phantom element of $H_i(F^e(M.))$.
\end{proof}

\begin{proof}[Proof (\ref{it:LMphAPX}).]
Let $[z] \in H_i(F^e(N.))$.  For any fixed $e' \geq e_0$, since $\delta_i^{(q',q'')}([z])$
is phantom for all $e'' \geq e_0$, it follows from Lemma~\ref{lem:deltaphzero} that
$\delta_i^{(q',q'')}([z]) = 0$ for all $e'' \geq e_0$.  Then since this holds for all $e' \geq
e_0$, the conclusion follows from part~(\ref{it:Mph})i.
\end{proof}

Finally in this section, we note the following characterization of stable phantom exactness.

\begin{prop}
A three-term complex $A \arrow{\alpha} B \arrow{\beta} C$ of finitely generated $R$-modules
is stably phantom exact if and only if for all $e$ (equiv. for all $e \gg 0$), the induced maps
$G^e(B / \im \alpha) \rightarrow G^e(C)$ are injective.  
\end{prop}

\begin{proof}
The first thing is to show that for any such 3-term complex, for every $e$ there is an
induced map $G^e(B / \im \alpha) \rightarrow G^e(C)$.  It is enough to show it for $e=0$, since
the tensor product of a complex with $^e R$ is still a complex.  The map will be induced from $\beta$, so
what needs to be shown is that $\beta((\im \alpha)^*) \subseteq 0^*_C$.  But if $y \in (\im \alpha)^*$,
then $c y^{q} \in \im F^e(\alpha)$ for all
$e \gg 0$, say $c y^q = \alpha(x_q)$, so $c \beta(y)^q = \beta(\alpha(x_q)) = 0$, which since this
holds for all $q \gg 0$, $\beta(y) \in 0^*_C$.  Thus, the induced maps are well-defined.

Now suppose that the complex is stably phantom exact, and take some $e$ and some $y \in F^e(B)$
such that $\beta(y) \in 0^*_{F^e(C)}$.  Then for all $q' \gg 0$, $\beta(c y^{q'}) = 0$, 
so $c y^{q'} \in (\im \alpha)^*$, and thus $c_1 y^{q' q_0} \in \im \alpha$, whence $y \in (\im \alpha)^*$.
This shows that the induced maps are injective for all $e$.

Conversely, suppose the induced maps are injective for $e \gg 0$.  Then take any $e$ and
any $y \in F^e(B)$ such that $\beta(y) = 0$.  Then for any $e' \gg 0$, the fact that
$\beta(y^{q'}) = \beta(y)^{q'} = 0$ implies by injectivity of the induced maps that
$y^{q'} \in (\im \alpha)^*_{F^{e+e'}(B)}$.  Thus $c y^{q' q_0} \in \im F^{e+e'+e_0}(\alpha)
= (\im F^e(\alpha))^{[q' q_0]}_{F^e(B)}$, so $y \in (\im F^e(\alpha))^*_{F^e(B)}$, showing that
the complex is stably phantom exact.
\end{proof}

\section{Analogues of right exactness, and tensoring with phantom regular sequences}
\label{sec:rtexact}
We show here that the
notion of right exactness generalizes well and that tensoring a short stably phantom
exact sequence with a sequence which is phantom regular on the third entry yields
a short stably phantom exact sequence.  See \cite[Proposition 1.1.4]{BH} for a non-phantom analogue.

\begin{lemma}\label{lem:rtphexact}
Let $R$ be a Noetherian ring of prime characteristic $p>0$.
Let $A \arrow{\alpha} B \arrow{\beta} C \rightarrow 0$ be a right stably phantom exact sequence
of finitely generated $R$-modules, and let $N$ be another finitely generated $R$-module.
Then \[
A \otimes_R N \arrow{\alpha \otimes N} B \otimes_R N \arrow{\beta \otimes N} C \otimes_R N \rightarrow 0
\]
is also a right stably phantom exact sequence of $R$-modules.
\end{lemma}
\begin{proof}
Let $R^t \arrow{f} R^s \rightarrow N \rightarrow 0$ be a finite free presentation of $N$.
Say $f(e_j) = \sum_{i=1}^s f_{i j} e_i$, $1 \leq j \leq t$, where $f_{i j} \in R$.  For any finitely
generated $R$-module $D$, let $g: D^t \rightarrow D^s$ be defined by setting
$g(d e_j) = \sum_{i=1}^s f_{i,j} d e_i$ for any $d \in D$, $1 \leq j \leq t$. 
Then by the right exactness of
tensor product, we can identify $D \otimes N$ with $\coker (1_D \otimes f) =
\coker g$.  Let $(\overline{b_1, \dotsc, b_s}) \in B \otimes N$
such that $(\overline{\beta(b_1), \dotsc, \beta(b_s)}) = 0 \in C \otimes N$.
That is, $b_i \in B$ for each $1 \leq i \leq s$, and \[
 \sum_{i=1}^s \beta(b_i) e_i = \tilde{f}_C(\sum_{j=1}^t c_j e_j) = 
 \sum_{i=1}^s (\sum_{j=1}^t f_{i,j} c_j) e_i
\] for some $c_1, \dotsc, c_t \in C$.  But each
$c_j = \beta(x_j)$ for some choices of $x_j \in B$, since $\beta$ is surjective.
Hence, $\beta(b_i - \sum_{j=1}^t f_{i,j} x_j) = 0$, $1 \leq i \leq s$.  
But then there is some $c \in R^o$ such that for any fixed $e \gg 0$,
there exist $a_i \in F^e(A)$, $1 \leq i \leq s$, such that $F^e(\alpha)(a_i) = 
c b_i^q - \sum_{j=1}^t c f_{i,j}^q x_j^q$.  Hence, \begin{align*}
\sum_{i=1}^s c b_i^q e_i &= \sum_{i=1}^s \left( F^e(\alpha)(a_i) e_i + \left(\sum_{j=1}^t c 
f_{i,j}^q x_j^q\right) e_i \right) \\
&= \sum_{i=1}^s F^e(\alpha)(a_i) e_i + \sum_{j=1}^t \left(\sum_{i=1}^s f_{i,j}^q 
c x_j^q e_i \right) \\
&= \sum_{i=1}^s F^e(\alpha)(a_i) e_i + \sum_{j=1}^t F^e(1_B \otimes f)(c x_j^q e_j)
\end{align*}
Hence, $c (\overline{b_1, \dotsc, b_s})^q = F^e(\alpha \otimes 1_N)(\overline{a_1,
\dotsc, a_s}) \in (\im (\alpha \otimes 1_N))^{[q]}_{B \otimes N}$.  So we have 
$\ker (\beta \otimes 1_N) \subseteq (\im (\alpha \otimes 
1_N))^*_{B \otimes N}$.  Since this all holds for arbitrary Frobenius powers of the
original sequence, it follows that the tensored sequence is stably phantom exact
at $B \otimes N$.  It is stably phantom exact at $C \otimes N$ by right-exactness
of the tensor product.

\end{proof}

\begin{lemma}\label{lem:phregtens}
Let $(R,\m)$ be a Noetherian local ring of prime characteristic $p>0$ containing a 
weak test element, let $0 \rightarrow A \arrow{\alpha}
B \arrow{\beta} C \rightarrow 0$ be a short stably phantom exact sequence of
finitely generated $R$-modules, and let $\bfx = x_1, \dotsc, x_n$
be a phantom $C$-regular sequence in $\m$.   Then \[
0 \rightarrow A / \bfx A \arrow{\bar{\alpha}} B / \bfx B \arrow{\bar{\beta}} C / \bfx C
\rightarrow 0
\]
is also a short stably phantom exact sequence of $R$-modules.
\end{lemma}

\begin{proof}
After Lemma~\ref{lem:rtphexact}, all that remains
is to show stable
phantom exactness at $A / \bfx A$.  Moreover, it suffices to show this for the case $n=1$, so $\bfx = x_1 = x$.  Let
$L. = K.(x; A)$, $M. = K.(x; B)$, $N. = K.(x; C)$, $\alpha_1 = \alpha_0 = \alpha$,
and $\beta_1 = \beta_0 = \beta$ Then it is easy to see that
\[
0 \rightarrow L. \arrow{\alpha.} M. \arrow{\beta.} N. \rightarrow 0
\]
satisfies the conditions of the Proposition~\ref{prop:pseudolesAPX}.  Moreover, the hypothesis
means that $N.$ has stably phantom homology at 1.  Hence by Proposition~\ref{prop:pseudolesAPX}
(\ref{it:NphAPX}), for 
any $[a] \in H_0(F^e(L.))$ which $H_0(F^e(\alpha.))$ sends to a phantom element of
$H_0(F^e(M.))$, $[a]$ itself must be phantom.  However, unrolling definitions, this
means that if $a \in F^e(A)$ and $F^e(\alpha)(a) \in (x^q F^e(B))^*_{F^e(B)}$, then 
$a \in (x^q F^e(A))^*_{F^e(A)}$.  In particular, $\ker F^e(\bar{\alpha}) \subseteq 
0^*_{F^e(\bar{A})}$.
\end{proof}

\section{Diagram chasing with exponents and Frobenius closures}\label{sec:frobdiag}
When investigating properties of tight closure, it is always tempting
to ask about analogous properties for Frobenius closure.  This latter
closure operation has a simpler definition, is easier to check, and
is known to commute with localization.  However, it also lacks some
of the nice properties of tight closure, so in the end it is probably
of less interest.  In this section, we shall sketch a method for recovering analogues of the
results of Section~\ref{sec:pseudolesAPX} for Frobenius closure in 
a much more intuitive way than the methods available for tight closure.

Recall the definition of Frobenius closure of submodules:  If $L \subseteq M$, then \[
L^F_M := \{z \in M \mid \exists q=p^e \text{ such that } z^q \in L^{[q]}_M \}.
\]

\begin{definition}
Let $R$ be a Noetherian ring of characteristic $p>0$, let $M$ and $N$
be $R$-modules, and let $\alpha: M \rightarrow N$ be an $R$-module map.
Then let \[
F^\infty(M) = \lim_{\rightarrow} F^e(M),
\]
where the maps in the direct limit system send $z \in F^e(M)$ to
$z^p \in F^{e+1}(M)$.  Let \[
F^\infty(\alpha): F^\infty(M) \rightarrow F^\infty(N)
\]
send the image of an element $z \in F^e(M)$ in $F^\infty(M)$ to the 
image of $F^e(\alpha)(z) \in F^e(N)$ in $F^\infty(N)$.  It is easy to show
that not only is this well-defined, but it 
makes $F^\infty$ into an endofunctor on the category of $R$-modules.  If 
$(C., d.)$ is a complex of $R$-modules, then $(F^\infty(C)., F^\infty(d).)$
is defined by $F^\infty(C)_i = F^\infty(C_i)$ and $F^\infty(d)_i = F^\infty(d_i)$.
The fact that $F^\infty$ is a functor is enough to show that this new sequence
is itself a complex.

Any element of $F^\infty(M)$ can be represented as the image of an element
$z \in F^e(M)$ for some $e$.  Denote the corresponding element of 
$F^\infty(M)$ by $\<z\>_e$.  Note that we always have $\<z\>_e = \<z^p\>_{e+1}$.
\end{definition}

\begin{lemma}\label{lem:Fexact}
Let $(C., d.)$ be a complex of $R$-modules over a Noetherian ring
of characteristic $p>0$.  Fix some integer $i$.  Then $F^\infty(C).$
is exact at $i$ if and only if $\ker F^e(d_i) \subseteq (\im 
F^e(d_{i+1}))^F_{F^e(C_i)}$ for all $e \geq 0$.
\end{lemma}

\begin{proof}
First suppose that $F^\infty(C).$ is exact at $i$.  Let $z \in \ker F^e(d_i)$.
Then \[
F^\infty(d)_i(\<z\>_e) = \<F^e(d_i)(z)\>_e = \<0\>_e = 0.
\]
Thus, by exactness, $\<z\>_e \in \im F^\infty(d)_{i+1}$.  Say 
$\<z\>_e = F^\infty(d)_{i+1}(\<y\>_n)$.  Without loss of generality $n \geq e$,
so $\<z^{p^{n-e}} - F^n(d_{i+1})(y)\>_n = 0$.  Then for some $t \geq n$, \[
z^{p^{t-e}} - F^t(d_{i+1})(y^{p^{t-n}}) = \left(z^{p^{n-e}} - F^n(d_{i+1})(y)
\right)^{p^{t-n}} = 0.
\]
That is, $z^{p^{t-e}} \in \im F^t(d_{i+1}) = (\im F^e(d_{i+1}))^{[p^{t-e}]}_{F^e(C_i)}$,
so that $z \in (\im F^e(d_{i+1}))^F_{F^e(C_i)}.$

Conversely, suppose that $\ker F^e(d_i) \subseteq (\im F^e(d_{i+1}))^F_{F^e(C_i)}$ for
all $e$, and let $\<z\>_n \in \ker F^\infty(d)_i$.  Then $\<F^n(d_i)(z)\>_n = 0$,
which means that for some $m \geq 0$, $F^{n+m}(d_i)(z^{p^m}) = 0$.  So we have that
\[
z^{p^m} \in \ker F^{n+m}(d_i) \subseteq (\im F^{n+m}(d_{i+1}))^F_{F^{n+m}(C_i)}.
\]
Thus, for some $t \geq 0$, $z^{p^{m+t}} \in (\im F^{n+m}(d_{i+1}))^{[p^t]}_{F^{n+m}(C_i)}
= \im F^{n+m+t}(d_{i+1}).$  Say $z^{p^{m+t}} = F^{n+m+t}(d_{i+1})(y).$  Then \[
\<z\>_n = \<z^{p^{m+t}}\>_{n+m+t} = \<F^{n+m+t}(d_{i+1})(y)\>_{n+m+t} = 
F^\infty(d)_{i+1}(\<y\>_{n+m+t}),
\]
which proves exactness of $F^\infty(C).$ at $i$.
\end{proof}

\begin{cor}\label{cor:Fexact}
Let $
0 \rightarrow L \stackrel{\alpha}{\rightarrow} M 
\stackrel{\beta}{\rightarrow} N \rightarrow 0
$
be a complex of modules over a local Noetherian ring $R$ of characteristic
$p>0$.  Then the following are equivalent:
\begin{enumerate}[(a)]
\item\label{it:a} The following three conditions hold:

\begin{enumerate}[i.]
\item\label{it:bsur} $\beta$ is surjective,
\item\label{it:kerim} $\ker F^e(\beta) \subseteq (\im F^e(\alpha))^F_{F^e(M)}$ for all $e$, and
\item\label{it:ainj} $\ker F^e(\alpha) \subseteq 0^F_{F^e(L)}$ for all $e$.
\end{enumerate}
\item\label{it:b} The sequence
$0 \rightarrow F^\infty(L) \stackrel{F^\infty(\alpha)}{\rightarrow} F^\infty(M) 
\stackrel{F^\infty(\beta)}{\rightarrow} F^\infty(N) \rightarrow 0$
is exact.
\end{enumerate}
\end{cor}
\begin{proof}
This is a special case of Lemma~\ref{lem:Fexact}.  The only part that remains
to be seen is the relationship of the surjectivity of $\beta$ to exactness
at $F^\infty(N)$.  If $\beta$ is surjective, then so is $F^e(\beta)$ for all
$e$ because of the right-exactness of Frobenius, from which it easily follows
that $F^\infty(\beta)$ is surjective.  On the other hand, if $F^\infty(\beta)$ is
surjective, then by Lemma~\ref{lem:Fexact}, $N \subseteq (\im \beta)^F_N$.
Hence $N \subseteq (\im \beta)^*_N$, which implies that $N = \im \beta$.
\end{proof}

From these results, along with the fact that direct limits are exact and thus
commute with homology, one can show that analogues of 
Lemma~\ref{lem:pseudoles1} and Proposition~\ref{prop:pseudolesAPX} hold for Frobenius closure.

\noindent \textbf{Remark:} The same kind of thing can be done with ``plus closure'', 
and it's even easier.  Recall: if $R$ is an integral domain and $j: L \hookrightarrow M$ is
an inclusion of $R$-modules,
then by definition the \emph{plus closure} of $L$ in $M$, denoted $L^+_M$, is
the set of all elements $x \in M$ such that for some module-finite extension $S$ of $R$,
$1 \otimes x \in L^S$, where $L^S$ is the image of the map $S \otimes j: S \otimes_R L 
\rightarrow S \otimes_R M$.  Equivalently, if we let $R^+$ denote a direct limit of
\emph{all} module-finite extensions of $R$ (equivalently, $R^+$ is the integral closure of $R$ in an 
algebraic closure of its fraction field), then $L^+_M$ is the preimage in $M$ of the image
of the map $R^+ \otimes j$ in $R^+ \otimes_R M$.  For a complex $C.$ of $R$-modules,
$C^+.$ is defined in the obvious way.

It should be clear to the reader that for a complex $(C., d.)$ and an integer $i$, $C^+.$ is exact at $i$ if and only if
$\ker (S \otimes_R d_i) \subseteq (\im (S \otimes_R d_{i+1}))^+_{S \otimes_R C_i}$ for
all module-finite extensions $S$ of $R$, which gives us a result analogous to Corollary~\ref{cor:Fexact}
and hence very easily versions of Lemma~\ref{lem:pseudoles1} and Proposition~\ref{prop:pseudolesAPX} for plus closure.

One would hope that there is a similarly straightforward
way to prove Lemma~\ref{lem:pseudoles1} and Proposition~\ref{prop:pseudolesAPX} as well.

\section*{Acknowledgements}
I warmly thank my graduate advisor, Professor Craig Huneke, for his support, guidance,
and ideas.  I also offer my thanks to the referee, who offered several helpful corrections and suggestions.

\providecommand{\bysame}{\leavevmode\hbox to3em{\hrulefill}\thinspace}
\providecommand{\MR}{\relax\ifhmode\unskip\space\fi MR }
\providecommand{\MRhref}[2]{%
  \href{http://www.ams.org/mathscinet-getitem?mr=#1}{#2}
}
\providecommand{\href}[2]{#2}

\end{document}